\documentclass[10pt]{article}

\usepackage{bm}
\usepackage{amsmath}
\usepackage{amsthm}
\usepackage{amsfonts}
\usepackage{amssymb}
\usepackage{amscd}
\usepackage{hyperref}
\usepackage{pstricks}
\usepackage{latexsym}
\usepackage{enumerate}
\usepackage{epsfig}
\usepackage{subfig}
\usepackage{graphicx}
\usepackage{caption}
\usepackage{framed,fancybox}
\usepackage{bbm}
\usepackage{multirow}
\usepackage[T1]{fontenc}
\usepackage{threeparttable}
\usepackage{rotating}
\usepackage{color,soul,xcolor}
\usepackage{physics}
\usepackage[version=4]{mhchem}
\usepackage{siunitx}
\usepackage{footnpag}			      	
\usepackage{longtable,tabularx}
\newcommand{\m}[1]{{\bf{#1}}}
\newcommand{\g}[1]{\boldsymbol #1}
\newcommand{\bb}[1]{\mathbb #1}
\newcommand{\C}[1]{{\cal {#1}}}
\newcommand{\mbb}[1]{\mathbb{#1}}
\newcommand{\R}{\mbb{R}}
\newcommand{\T}{^{\sf T}}    

\makeatletter
\newcommand{\thickhline}{%
    \noalign {\ifnum 0=`}\fi \hrule height 1pt
    \futurelet \reserved@a \@xhline
}
\newcolumntype{"}{@{\hskip\tabcolsep\vrule width 1pt\hskip\tabcolsep}}
\makeatother

\definecolor{shadecolor}{gray}{0.99}
{\endMakeFramed}

\newenvironment{shadedframe}{%
 \MakeFramed {\FrameRestore}}
{\endMakeFramed}

\long\def\symbolfootnote[#1]#2{\begingroup\def\thefootnote{\fnsymbol{footnote}}\footnote[#1]{#2}\endgroup}
\def\qed{\hfill{$\vcenter{\hrule height1pt \hbox{\vrule width1pt height5pt
    \kern5pt \vrule width1pt} \hrule height1pt}$} \medskip}
\title{\bf Method for Solving Bang-Bang and Singular Optimal \\ Control Problems using Adaptive Radau Collocation}

\author{Elisha R.~Pager\thanks{Ph.D.~Candidate, Department of Mechanical and Aerospace Engineering, University of Florida, Gainesville, Florida 32611-6250. Email: epager@ufl.edu.} \\
Anil V.~Rao\thanks{Professor, Department of Mechanical and Aerospace Engineering, University of Florida, Gainesville, FL 32611-6250. E-mail:  anilvrao@ufl.edu.  Corresponding Author.} \vspace{12pt} \\ {\em University of Florida} \\ {\em Gainesville, FL 32611}
}

\date{}

\begin{document}
\maketitle
\thispagestyle{empty}

\begin{abstract}
  A method is developed for solving bang-bang and singular optimal control problems using adaptive Legendre-Gauss-Radau (LGR) collocation.  The method is divided into several parts.  First, a structure detection method is developed that identifies switch times in the control and analyzes the corresponding switching function for segments where the solution is either bang-bang or singular.  Second, after the structure has been detected, the domain is decomposed into multiple domains such that the multiple-domain formulation includes additional decision variables that represent the switch times in the optimal control.  In domains classified as bang-bang, the control is set to either its upper or lower limit.  In domains identified as singular, the objective function is augmented with a regularization term to avoid the singular arc.  An iterative procedure is then developed for singular domains to obtain a control that lies in close proximity to the singular control.  The method is demonstrated on four examples, three of which have either a bang-bang and/or singular optimal control while the fourth has a smooth and nonsingular optimal control.  The results demonstrate that the method of this paper provides accurate solutions to problems whose solutions are either bang-bang or singular when compared against previously developed mesh refinement methods that are not tailored for solving nonsmooth and/or singular optimal control problems, and produces results that are equivalent to those obtained using previously developed mesh refinement methods for optimal control problems whose solutions are smooth. \\
\end{abstract}


\renewcommand{\baselinestretch}{1}
\normalsize\normalfont

\renewcommand{\baselinestretch}{2}
\normalsize\normalfont 

\section{Introduction}

Optimal control problems arise in many engineering applications due to the need to optimize the performance of a controlled dynamical system.  In general, optimal control problems do not have analytic
solutions and must be solved numerically.  A key challenge in solving an optimal control problem numerically arises due to the fact that most optimal control problems are subject to constraints on the system. These constraints often take the form of path constraints where limits are imposed on functions of either the control and/or the state.  Constrained optimal control problems often have nonsmooth solutions, where the nonsmoothness arises in the forms of instantaneous switches in the control or switches between activity and inactivity in the state path constraints.  Moreover, many constrained optimal control problems have solutions that lie on one or more singular arcs.  The existence of a singular arc makes solving constrained optimal control problems even more challenging because Pontryagin's minimum principle (that is, the first and second-order optimality conditions) fail to yield a complete solution along the singular arc. As a result, when applying a computational method to a problem whose solution lies on a singular arc, standard methods produce nonsensical results. This research is motivated by the importance of solving optimal control problems whose solutions are nonsmooth and singular.

Numerical methods for optimal control fall into two broad categories: indirect methods and direct methods. In an indirect method, the first-order variational optimality conditions are derived, and the optimal control problem is converted to a Hamiltonian boundary-value problem (HBVP). The HBVP is then solved numerically using a differential-algebraic equation solver. In a direct method, the state and control are approximated, and the optimal control problem is transcribed into a finite-dimensional nonlinear programming problem (NLP)~\cite{Betts2010}. The NLP is then solved numerically using
well-developed software such as {\em SNOPT}~\cite{Gill2002} or {\em IPOPT}~\cite{Biegler2008}.

Over the past two decades, a particular class of direct methods, called direct collocation methods, has been used extensively for solving continuous optimal control problems. A direct collocation method is
an implicit simulation method where the state and control are parameterized, and the constraints in the continuous optimal control problem are enforced at a specially chosen set of collocation points.
In more recent years, a great deal of research has been carried out in the area of {\em Gaussian quadrature orthogonal collocation methods}~\cite{Benson2006,Garg2010,Garg2011a,Garg2011b}.  In a Gaussian quadrature collocation method, the state is approximated using a basis of Lagrange polynomials where the support points of the polynomials are chosen to be the points associated with a Gaussian quadrature.   The most well developed Gaussian quadrature methods employ either Legendre-Gauss (LG) points~\cite{Benson2006,Rao2010}, Legendre-Gauss-Radau (LGR) points~\cite{Garg2010,Garg2011a,Garg2011b,Kameswaran2008,Patterson2015}, or Legendre-Gauss-Lobatto (LGL) points~\cite{Elnagar1995}.  More recently, a convergence theory has been developed to show that under certain assumptions of smoothness and coercivity, an $hp$--Gaussian quadrature method employing either LG or LGR collocation points will converge to a local minimizer of the optimal control problem at an exponential rate~\cite{Hager2016a,Hager2016b,Hager:2018,Hager:2019}. For these reasons, direct methods should be explored in the effort to develop novel computational methods for singular optimal control.

Computational issues arise when a solution to an optimal control problem is either nonsmooth or singular. The difficulty with such optimal control problems is twofold.  First, the precise locations of any discontinuities and the structure of the control must be identified.  Static mesh refinement methods that employ Gaussian quadrature have been developed recently as an initial attempt to locate (approximately) discontinuities in the solution~\cite{Liu2015,Gong2008,Miller2021}.   Even more recently, the idea of using variable mesh refinement methods and structure detection methods has been developed~\cite{Schlegel2004,Schlegel2006,Wang2014,Chen2016,Chen2019}. Unlike static mesh refinement methods such as those found in Refs.~\cite{Patterson2015,Darby2010,Liu2018}, variable mesh refinement methods work by including parameters in the optimization that define the location of the discontinuities. The methods in~Refs.~\cite{Schlegel2004,Schlegel2006,Wang2014} use the Lagrange multipliers to detect the switch point locations in the control structure and then place variable mesh points to represent the switch times in the NLP, while~Refs.~\cite{Chen2016,Chen2019} use the switching function and a sensitivity analysis to place moving finite elements at the switch point locations. Furthermore, Ref.~\cite{Agamawi2020} describes a mesh refinement method for solving bang-bang optimal control problems based on the switching function associated with the Hamiltonian.  More recently, in Ref.~\cite{Aghaee2021} a switch point algorithm was developed for optimizing over the locations of switch points in a nonsmooth control solution, but a priori knowledge of the switch points existence is required.  Finally, methods that utilize structure detection on a static mesh are described in Refs.~\cite{Kaya2003,Mehrpouya2021}.

The second difficulty in solving optimal control problems with nonsmooth or singular solutions arises when the optimal control is singular.  Several approaches have been developed for solving singular optimal control problems using both indirect and direct methods. A majority of these methods typically employ either a regularization approach or use of the optimality conditions with an indirect method to solve for the singular control (see Refs.~\cite{SoledadAronna2013,Mehra1972}).  A regularization method transforms the singular control problem into a series of nonsingular problems by minimizing the sum of the original objective and a regularization term, where the regularization term is a quadratic function of the control.  Regularization approaches have been implemented using dynamic programming, indirect methods, direct methods, and nested indirect/direct approaches as described in Refs.~\cite{Jacobson1970,Maurer1976,AndrsMartnez2018,Caponigro2018}.  More recent implementations of regularization based techniques include the uniform trigonometrization method (UTM) developed in~\cite{Mall2020}, the use of a continuation method of a regularized term in~\cite{Fabien2013}, and the total variation based regularization approach in~\cite{Aghaee2021}.  Aside from regularization based approaches, research has also been conducted on the use of low-order representations of the control including straight line, monotonic, and nonmonotonic function approximations to reduce the oscillations and numerical challenges observed with singular arcs~\cite{AndrsMartnez2020,Maga1997}.

Motivated by the prevalence of bang-bang and singular arcs in optimal control solutions and the need for a general method that can handle such problems, this paper describes a new method for detecting and solving optimal control problems whose solutions are bang-bang and singular.  The method described in this paper consists of the following parts.  First, a multiple-domain reformulation of the LGR collocation method is developed that enables partitioning the problem into segments that are categorized as either regular, bang-bang, or singular. This multiple-domain LGR collocation reformulation partitions the time horizon into domains such that additional decision variables are introduced that correspond to the endpoints of each newly created domain.  Second, a method is developed that detects the structure of the optimal control (that is, the method detects segments where the solution is either regular, bang-bang, or singular).  Third, a regularization method (inspired by the approach used in Ref.~\cite{Jacobson1970}) is employed to solve for the control in any domain where the control is categorized as singular.  The method presented in~\cite{Jacobson1970} implements a regularization procedure over the entire time domain using differential dynamic programming, whereas in this current work a similar regularization procedure is implemented only in the intervals denoted as singular. This difference in implementation results in a more accurate solution of the control. It is noted that, because bang-bang and singular optimal control problems frequently have Hamiltonians that are affine in the control, the structure detection method developed in this paper is designed to work on those components of the control that appear linearly in the Hamiltonian.  Moreover, switch times in the control (which then lead to the partitioning into multiple domains) are identified using the jump function mesh refinement method described in Ref.~\cite{Miller2021}.  Using the aforementioned structure detection method, the control in the newly created domains is then either left free (in the case of a regular domain), set equal to one of its limits (in the case of a bang-bang domain), or is determined via the aforementioned regularization method (in the case of a singular domain).  


The contributions of this work are as follows. First, using jump function approximations provide an accurate way to determine the number of discontinuities along with accurate estimates of the locations of these discontinuities.  Second, using the switching function and the Hamiltonian enables determining those intervals where the control is either bang-bang or singular. Third, the method automatically partitions the solution into domains based on the results of the structure detection method.  Fourth, the method does not require any a priori knowledge of the structure in the optimal control or whether the optimal control is bang-bang or singular.  Fifth, in this paper the use of regularization methods is extended to direct collocation methods.  In particular, the multiple-domain partition of the solution obtained from the structure detection method enables regularizing only over those domains where the control is singular.  Consequently, within a singular domain the regularization leads to a control that lies in close proximity to the singular control while simultaneously eliminating the need to derive the singular control conditions (where deriving such conditions may prove to be intractable depending upon the problem).  The performance of the method developed in this paper is demonstrated on four examples.  The optimal control for each of the first three of these examples is either bang-bang and/or singular, while the optimal control for the fourth example is smooth.  The numerical results obtained of the first three examples demonstrate that the method of this paper produces significantly more accurate results when compared against mesh refinement methods that are not developed for solving optimal control problems whose solutions are nonsmooth or singular.  Finally, the numerical results of the fourth example demonstrate the method of this paper correctly identifies when a solution is smooth and applies only static mesh refinement in order to obtain a solution.  As a result, results comparable to those obtained using previously developed mesh refinement methods are obtained.

The remainder of the paper is organized as follows. Section~\ref{sect:probForm} introduces the Bolza optimal control problem and the necessary conditions for optimality. Section~\ref{sect:LGR} describes the multiple-domain Legendre-Gauss-Radau collocation used to transcribe the multiple-domain Bolza optimal control problem. Section~\ref{sect:SOCP} provides a brief overview of optimal control problems whose solutions are nonsmooth. Section~\ref{sect:method} details the method for solving bang-bang and singular optimal control problems.  Section~\ref{sect:examples} provides numerical solutions obtained by demonstrating the method on four examples. 
Section \ref{sect:limitations} describes limitations of the method.  Finally, Section~\ref{sect:conc} provides conclusions on this research.  

\section{Bolza Optimal Control Problem\label{sect:probForm}}

Without loss of generality, consider the following single-phase optimal control problem in Bolza form defined on the time horizon $t\in[t_0,t_f]$. Determine the state $\m{x}(t)\in\bb{R}^{n_x}$, the control $\m{u}(t)\in\bb{R}^{n_u}$, and the terminal time $t_f\in\bb{R}$ that minimize the objective functional
\begin{equation}\label{eq:single-cost-t}
  		\C{J}=\C{M}(\m{x}(t_0),t_0,\m{x}(t_f),t_f)+ \int_{t_0}^{t_f} \C{L}(\m{x}(t),\m{u}(t), t)~dt, 
\end{equation}
subject to the dynamic constraints
\begin{equation}\label{eq:single-dyn-t}
  \frac{d\m{x}(t)}{dt} \equiv \m{\dot{x}}(t) = \m{a}(\m{x}(t),\m{u}(t),t),
\end{equation}
the control constraints
\begin{equation}\label{eq:single-control-t}
	\m{u}_{\min} \leq \m{u}(t) \leq \m{u}_{\max},
\end{equation}
and the boundary conditions
\begin{equation}\label{eq:single-bc-t}
  \m{b}(\m{x}(t_0),t_0,\m{x}(t_f),t_f) = \textbf{0},
\end{equation}
where the functions $\C{M}$, $\C{L}$, $\m{a}$, $\m{b}$, and $\m{c}$ are defined by the mappings
\begin{displaymath}
	\begin{array}{lcl}
		\C{M} & : & \R^{n_x} \times \R \times \R^{n_x} \times \R \rightarrow \R, \\
		\C{L} & : & \R^{n_x}  \times \R^{n_u} \times \R \rightarrow \R, \\
		\m{a} & : & \R^{n_x} \times \R^{n_u} \times \R \rightarrow \R^{n_x}, \\
		\m{b} & : & \R^{n_x} \times \R \times \R^{n_x} \times \R \rightarrow \R^{n_b}. \\
	\end{array}
\end{displaymath}

The Bolza optimal control problem given in Eqs.~\eqref{eq:single-cost-t}--\eqref{eq:single-bc-t} gives rise to the following first-order calculus of variations~\cite{Athans2013,Kirk2004,Bryson1975} conditions:
\begin{eqnarray}
  \dot{\m{x}}(t) & = & \phantom{-}\left[\frac{\partial \C{H}}{\partial \g{\lambda}}\right]\T = \phantom{-}\C{H}_{\g{\lambda}}\T, \label{eq:Hamiltonian-system1} \\
  \dot{\g{\lambda}}(t) & = & - \left[\frac{\partial \C{H}}{\partial \m{x}}\right]\T = -\C{H}_{\m{x}}\T, \label{eq:Hamiltonian-system2} \\
	\m{0} &=& \phantom{-}\frac{\partial \C{H} }{\partial \m{u}} = \phantom{-}\C{H}_{\m{u}}, \label{eq:pontryagin-strong-form}
\end{eqnarray}
where $\g{\lambda}(t)\in\R^{n_x}$ is the costate, 
\begin{equation}\label{eq:controlHamiltonian}
\C{H}(\m{x}(t),\m{u}(t),\g{\lambda}(t),t) = \C{L}(\m{x}(t),\m{u}(t),t)  + \g{\lambda}\T(t)\m{a}(\m{x}(t),\m{u}(t),t),
\end{equation}
is the augmented Hamiltonian, 
and $\C{U}$ is the admissible control set.  Finally, the transversality conditions are given by
\begin{equation}\label{eq:transversality-conditions}
\begin{array}{lclclcl}
  \g{\lambda}(t_0) & = & -\displaystyle\frac{\partial{\C{M}}}{\partial{\m{x}(t_0)}}  + \g{\nu}\T
  \displaystyle\frac{\partial{\m{b}}}{\partial{\m{x}(t_0)}}  & , & 
  \g{\lambda}(t_f) & = & \displaystyle\frac{\partial{\C{M}}}{\partial{\m{x}(t_f)}} - \g{\nu}\T \displaystyle\frac{\partial{\m{b}}}{\partial{\m{x}(t_f)}}
\end{array}
\end{equation}
\begin{equation}\label{eq:initial-and-final-H}
\begin{array}{lclclcl}
  \C{H}(t_0) & = & \displaystyle\frac{\partial{\C{M}}}{\partial{t_0}} - \g{\nu}\T\frac{\partial{\m{b}}}{\partial{t_0}} & , & \C{H}(t_f) &  = & -\displaystyle\frac{\partial{\C{M}}}{\partial{t_f}}  + \g{\nu}\T \frac{\partial{\m{b}}}{\partial{t_f}} 
\end{array}
\end{equation}
where $\g{\nu}$ is the Lagrange multiplier associated with the boundary conditions.  Equations \eqref{eq:Hamiltonian-system1} and \eqref{eq:Hamiltonian-system2} form what is classically known as a {\em Hamiltonian system}~\cite{Athans2013,Kirk2004}.  The conditions in Eqs.~\eqref{eq:transversality-conditions} and~\eqref{eq:initial-and-final-H} are called {\em transversality conditions}~\cite{Athans2013,Kirk2004,Bryson1975} on the boundary values of the costate,
For some problems, the control cannot be uniquely determined, either implicitly or explicitly, from the optimality conditions given in Eqs.~\eqref{eq:Hamiltonian-system1}--\eqref{eq:initial-and-final-H}. In such cases, the weak form of Pontryagin's minimum principle can be used which solves for the permissible control that minimizes the Hamiltonian in Eq.~\eqref{eq:controlHamiltonian}. If $\m{u} \in [u_{\min},u_{\max}]$ is the set of permissible controls, then Pontryagin's minimum principle states that the optimal control, $\m{u}$, satisfies the condition
\begin{equation}\label{eq:PMP}
  \C{H}(\m{x}^*(t),\m{u}^*,\g{\lambda}^*(t),t) \leq \C{H}(\m{x}^*(t),\m{u},\g{\lambda}^*(t),t), \quad \m{u} \in [u_{\min},u_{\max}].
\end{equation}
The Hamiltonian system, together with the original boundary conditions and the costate transversality conditions, forms a {\em Hamiltonian boundary-value problem}  (HBVP)~\cite{Athans2013,Kirk2004,Bryson1975}.  Any solution $(\m{x}^*(t),\m{u}^*,\g{\lambda}^* (t),\g{\nu}^*)$ to the HBVP is called an {\em extremal} solution. 


\section{Multiple-Domain Legendre-Gauss-Radau Collocation\label{sect:LGR}}

In this paper, the previously developed $hp$--adaptive Legendre-Gauss-Radau (LGR) collocation method \cite{Garg2010,Garg2011a,Garg2011b,Kameswaran2008,Patterson2015} is used to approximate the optimal control problem (where the term LGR collocation will be used from this point onwards to mean $hp$--adaptive Legendre-Gauss-Radau collocation).  LGR collocation is used because it has been shown to converge at an exponential rate to a local solution of the optimal control problem for problems where the solution is smooth \cite{Hager2016a,Hager2016b,Hager:2018,Hager:2019}. The focus of this paper, however, is on solving optimal control problems whose solutions are nonsmooth and/or singular.  As a result, modifications to the standard LGR formulation are made.  Specifically, a multiple-domain reformulation of LGR collocation is developed as described in the remainder of this section.  

The multiple-domain formulation of LGR collocation divides the domain $t\in[t_0,t_f]$ into distinct partitions such that the endpoints of each partition are decision variables.  The division into domains is obtained using a structure decomposition method as described in Section \ref{sect:sructureDet}. The continuous-time Bolza optimal control problem described in Eqs.~\eqref{eq:single-cost-t}--\eqref{eq:single-bc-t} is discretized using collocation at the Legendre-Gauss-Radau (LGR) points~\cite{Garg2010,Garg2011a,Garg2011b,Kameswaran2008}.
The time horizon $t\in[t_0,t_f]$ may be divided into $D$ time domains, $\C{P}_d=[t_s^{[d-1]},t_s^{[d]}]\subseteq[t_0,t_f],~d\in\{1,\ldots,D\}$, such that
\begin{equation}\label{eq:time-domain-properties}
\bigcup_{d=1}^{D}\C{P}_d=[t_0,t_f] ~, \quad\bigcap_{d=1}^{D}\C{P}_d=\{t_s^{[1]},\ldots,t_s^{[D-1]}\}~,
\end{equation}
where $t_s^{[d]},~d\in\{1,\ldots,D-1\}$ are the domain interface variables of the problem, $t_s^{[0]}=t_0$, and $t_s^{[D]} = t_f$. Thus, in the case where $D=1$ the phase consists of only a single domain $\C{P}_1=[t_0,t_f]$ and $\{t_s^{[1]},\ldots,t_s^{[D-1]}\}=\emptyset$.  
\begin{equation}\label{eq:affine-transformation}
  \begin{array}{lcl}
    t & = & \displaystyle \frac{t_s^{[d]}-t_s^{[d-1]}}{2}\tau + \frac{t_s^{[d]}+t_s^{[d-1]}}{2}~, \vspace{3pt} \\
    \tau & = &  \displaystyle 2\frac{t-t_s^{[d-1]}}{t_s^{[d]}-t_s^{[d-1]}}-1~.
  \end{array}
\end{equation}
The interval $\tau\in[-1,+1]$ for each domain $\C{P}_d$ is then divided into $K$ mesh intervals, $\C{I}_k=[T_{k-1},T_k]\subseteq [-1,+1],\;k\in\{1,\ldots,K\}$ such that
\begin{equation}\label{eq:mesh-interval-properties}
\bigcup_{k=1}^K\C{I}_k=[-1,+1] ~, \quad\bigcap_{k=1}^K\C{I}_k=\{T_1,\ldots,T_{K-1}\}~,
\end{equation}
and $-1=T_0<T_1<\ldots<T_{K-1}<T_K=+1$.  For each mesh interval, the LGR points used for collocation are defined in the domain of $[T_{k-1},T_k]$ for $k\in\{1,\ldots,K\}$.  The state of the continuous optimal control problem is then approximated in mesh interval $\C{I}_k,\;k\in\{1,\ldots,K\}$, as 
\begin{equation}\label{eq:LGR-state-approximation}
\m{x}^{(k)}(\tau)  \approx \m{X}^{(k)}(\tau) = \sum_{j=1}^{N_k+1} \m{X}_{j}^{(k)}
\ell_{j}^{(k)}(\tau)~, \quad  \ell_{j}^{(k)}(\tau) = \prod_{\stackrel{l=1}{l\neq j}}^{N_k+1}\frac{\tau-\tau_{l}^{(k)}}{\tau_{j}^{(k)}-\tau_{l}^{(k)}}~, 
\end{equation}  
where
$\ell_{j}^{(k)}(\tau)$ for $ j\in\{1,\ldots,N_k+1\}$ is a basis of Lagrange polynomials on $\C{I}_k$, $\left(\tau_1^{(k)},\ldots,\tau_{N_k}^{(k)}\right)$ are the set of $N_k$ Legendre-Gauss-Radau (LGR) collocation points in the interval $[T_{k-1},T_k)$, $\tau_{N_k+1}^{(k)}=T_k$ is a non-collocated support point, and $\m{X}_{j}^{(k)} \equiv  \m{X}^{(k)}(\tau_j^{(k)})$.  Differentiating $\m{X}^{(k)}(\tau)$ in Eq.~\eqref{eq:LGR-state-approximation} with respect to $\tau$ gives
\begin{equation}\label{eq:LGR-diff-state-approximation}
  \frac{d\m{X}^{(k)}(\tau)}{d\tau} = \sum_{j=1}^{N_k+1}\m{X}_{j}^{(k)}\frac{d\ell_j^{(k)}(\tau)}{d\tau}~.
\end{equation}
The dynamics are then approximated at the $N_k$ LGR points in mesh
interval $k\in\{1,\ldots,K\}$ as
\begin{equation}\label{eq:LGR-defect}
 \sum_{j=1}^{N_k+1}D_{lj}^{(k)} \m{X}_j^{(k)} - \frac{t_f-t_0}{2}\m{a}\left(\m{X}_l^{(k)},\m{U}_l^{(k)},t (\tau_l^{(k)},t_0,t_f)\right) = \m{0} ~,\quad l \in \{1,\ldots,N_k\}~,
\end{equation}
where 
\begin{equation*}
  D_{lj}^{(k)} = \frac{d\ell_j^{(k)}(\tau_l^{(k)})}{d\tau}~,\quad l \in \{1,\ldots,N_k\}~,~ j \in \{1,\ldots,N_k+1\}~,
\end{equation*}
are the elements of the $N_k\times (N_k+1)$ {\em Legendre-Gauss-Radau differentiation matrix} in mesh interval $\C{I}_k$, $\;k\in\{1,\ldots,K\}$, and $\m{U}_l^{(k)}$ is the approximation of the control at the $l^{th}$ collocation point in mesh interval $\C{I}_k$. The time variables $t_0$ and $t_f$ in Eq.~\eqref{eq:LGR-defect} represent the initial and final domain interface variables, $t_s^{[d-1]}$ and $t_s^{[d]}$, on the domain $\C{P}_d$.  It is noted that continuity in the state and time between mesh intervals $\C{I}_{k-1}$ and $\C{I}_{k}$, $k\in\{1,\ldots,K\}$, is enforced by using the same variables to represent $\m{X}_{N_{k-1}+1}^{(k-1)}$ and $\m{X}_{1}^{(k)}$, while continuity in the state between the domains $\C{P}_{d-1}$ and $\C{P}_{d}$, $d\in\{2,\ldots,D\}$, is achieved by using the same variables to represent $\m{X}_{N^{[d-1]}+1}^{[d-1]}$ and $\m{X}_{1}^{[d]}$ where the superscript $[d]$ is used to denote the $d^{th}$ time domain, $\m{X}_{j}^{[d]}$ denotes the value of the state approximation at the $j^{th}$ discretization point in the time domain $\C{P}_d$, and $N^{[d]}$ is the total number of collocation points used in time domain $\C{P}_d$ computed by
\begin{equation}\label{eq:N-sum-domain}
N^{[d]} = \sum_{k=1}^{K^{[d]}} N_k^{[d]}~.
\end{equation}

The Legendre-Gauss-Radau approximation of the multiple-domain optimal control problem results in the following nonlinear programming problem (NLP).  Minimize the objective function
\begin{equation}\label{eq:NLP-cost}
  \C{J}=\C{M}(\m{X}_1^{[1]},t_0,\m{X}_{N^{[D]}+1}^{[D]},t_f)+ \sum_{d=1}^{D} \frac{t_s^{[d]}-t_s^{[d-1]}}{2}\left[\m{w}^{[d]}\right]\T\m{L}^{[d]}~,
\end{equation}
subject to the collocated dynamic constraints
\begin{equation}\label{eq:NLP-defect}
\g{\Delta}^{[d]} =
\m{D}^{[d]}\m{X}^{[d]} - \frac{t_s^{[d]}-t_s^{[d-1]}}{2}\m{A}^{[d]}
=\m{0}
~,\quad  d \in \{1,\ldots,D\}~,
\end{equation}
the control constraints
\begin{equation}\label{eq:NLP-control}
\m{u}_{\min} \leq
\m{U}_{j}^{[d]}
\leq \m{u}_{\max}
~,\quad  j \in \{1,\ldots,N^{[d]}\}~,~ d \in \{1,\ldots,D\}~,
\end{equation}
the boundary conditions
\begin{equation}\label{eq:NLP-bc}
  \m{b}(\m{X}_1^{[1]},t_0,\m{X}_{N^{[D]}+1}^{[D]},t_f) = \textbf{0}~,
\end{equation}
and the continuity constraints
\begin{equation}\label{eq:NLP-continuity}
\m{X}_{N^{[d-1]}+1}^{[d-1]} = \m{X}_{1}^{[d]}~,~ d \in \{2,\ldots,D\}~,
\end{equation}
noting that Eq.~\eqref{eq:NLP-continuity} is implicitly satisfied by employing the same variable in the NLP for $\m{X}_{N^{[d-1]}+1}^{[d-1]}$ and $\m{X}_{1}^{[d]}$.
The matrices in Eqs.~\eqref{eq:NLP-cost}--\eqref{eq:NLP-defect} are defined as follows
\begin{equation}\label{eq:NLP-A-def}
\m{A}^{[d]} = \begin{bmatrix}
\m{a}\left(\m{X}_{1}^{[d]},\m{U}_{1}^{[d]},t_1^{[d]}\right)\\
\vdots\\
\m{a}\left(\m{X}_{N^{[d]}}^{[d]},\m{U}_{N^{[d]}}^{[d]},t_{N^{[d]}}^{[d]}\right)
\end{bmatrix}
\in \mathbb{R}^{N^{[d]}~\times~n_x}~,
\end{equation}
\begin{equation}\label{eq:NLP-L-def}
\m{L}^{[d]} = \begin{bmatrix}
\C{L}\left(\m{X}_{1}^{[d]},\m{U}_{1}^{[d]},t_1^{[d]}\right)\\
\vdots\\
\C{L}\left(\m{X}_{N^{[d]}}^{[d]},\m{U}_{N^{[d]}}^{[d]},t_{N^{[d]}}^{[d]}\right)
\end{bmatrix}
\in \mathbb{R}^{N^{[d]}~\times~1}~,
\end{equation}
$\m{D}^{[d]} \in \mathbb{R}^{N^{[d]}~\times~[N^{[d]}+1]}$ is the LGR differentiation matrix in time domain $\C{P}_d,~d \in\{1,\ldots,D\}$, and $\m{w}^{[d]} \in \mathbb{R}^{N^{[d]}~\times~1}$ are the LGR weights at each node in time domain $\C{P}_d,~ d \in\{1,\ldots,D\}$.  It is noted that $\m{a} \in \mathbb{R}^{1~\times~n_x}$ and $\C{L} \in \mathbb{R}^{1~\times~1}$ correspond, respectively, to the vector fields that define the right-hand side of the dynamics and the integrand of the optimal control problem.  Additionally, the state matrix, $\m{X}^{[d]}\in \mathbb{R}^{[N^{[d]}+1]~\times~n_x}$, and the control matrix, $\m{U}^{[d]}\in \mathbb{R}^{N^{[d]}~\times~n_u}$, in time domain $\C{P}_d,~ d \in \{1,\ldots,D\}$, are formed as 
\begin{equation}\label{eq:NLP-Var-Mat}
\m{X}^{[d]} = \begin{bmatrix}
\m{X}^{[d]}_{1}\\
\vdots \\
\m{X}^{[d]}_{N^{[d]}+1}
\end{bmatrix}
\text{ and }
\m{U}^{[d]} = 
\begin{bmatrix}
\m{U}^{[d]}_{1}  \\
\vdots	\\
\m{U}^{[d]}_{N^{[d]}}
\end{bmatrix}~,
\end{equation}
respectively, where $n_u$ is the number of control components and $n_x$ is the number of state components in the problem.  

\subsection{Costate Estimation\label{sect:costates}}

Estimates of the costate may be obtained at each of the discretization points in the time domain $\C{P}_d, d\in\{1,\ldots,D\}$ using the transformation \cite{Garg2010,Garg2011a,Garg2011b},
\begin{equation}\label{eq:NLP-costates} 
	\begin{array}{lcl}
		\g{\lambda}^{[d]} &=& (\m{W}^{[d]})^{-1}\g{\Lambda}^{[d]}~,  \\
		\g{\lambda}_{N^{[d]}+1}^{[d]} &=& (\m{D}_{N^{[d]}+1}^{[d]})\T \g{\Lambda}^{[d]}~,
	\end{array}
\end{equation}
where $\g{\lambda}^{[d]}\in\mathbb{R}^{N^{[d]}~\times~n_x}$ is a matrix of the costate estimates at the collocation points in time domain $\C{P}_d$, $\m{W}^{[d]} = \text{diag}(\m{w}^{[d]})$ is a diagonal matrix of the LGR weights at the collocation points in time domain $\C{P}_d$, $\g{\Lambda}^{[d]}\in\mathbb{R}^{N^{[d]}~\times~n_x}$ is a matrix of the NLP multipliers obtained from the NLP solver corresponding to the defect constraints at the collocation points in time domain $\C{P}_d$ , $\g{\lambda}_{N^{[d]}+1}^{[d]}\in\mathbb{R}^{1~\times~n_x}$ is a row vector of the costate estimates at the non-collocated end point in time domain $\C{P}_d$, and $\m{D}_{N^{[d]}+1}^{[d]}\in\mathbb{R}^{N^{[d]}~\times~1}$ is the last column of the LGR differentiation matrix in time domain $\C{P}_d$.

The aforementioned multiple-domain LGR formulation is summarized as follows.  First, a single {\em phase} problem on $t \in [t_0, t_f]$ is divided into $D$ {\em domains}, $\C{P}_d = [t_s^{[d-1]},t_s^{[d]}] , \: d \in \{1,\ldots,D\}$.  Each of the $D$ domains are then mapped to the interval $\tau^{[d]}\in[-1,+1],~d\in\{1,\ldots,D\}$. The interval $\tau^{[d]}\in[-1,+1],~d\in\{1,\ldots,D\}$ for each domain is then divided into $K$ {\em mesh intervals}, $\C{I}_k=[T_{k-1},T_k]\subseteq [-1,+1],\;k\in\{1,\ldots,K\}$.  Finally, the intersection of each domain is determined by the {\em domain interface variables}, $t^{[d]}_s,~d\in\{1,\ldots,D-1\}$.


\section{Nonsmooth and Singular Optimal Control\label{sect:SOCP}}

The term nonsmooth is used to denote the optimal control as displaying both nonsmooth and singular behavior.  By definition a singular arc occurs when Eq.~\eqref{eq:PMP} fails to uniquely describe an optimal control; for example, the set of minimizers in Eq.~\eqref{eq:PMP} form an interval which contains the optimal control~\cite{Bryson1975}. In this case, knowledge of the interval does not uniquely describe the optimal control itself. This phenomena can occur in many situations but is most common when the dynamics are linear in the control and the control is bounded, or the Hamiltonian is not an explicit function of time. It should be noted that singular arcs can also occur in other situations, but in order to provide structure to the method developed in this paper, only problems that fall into the aforementioned categories will be considered.

For simplicity, in the discussion that follows it is assumed that the control is a scalar (that is, $u(t)\in\bb{R}$).  Note, however, that without loss of generality the discussion below can be extended to multiple control components, and the use of the method developed in this paper on problems with multiple control components is demonstrated in the examples provided in Section~\ref{sect:examples}. Suppose the optimal control problem described in Eqs.~\eqref{eq:single-cost-t}--\eqref{eq:single-bc-t} is nonsmooth and singular as defined by the assumptions mentioned previously. 
The dynamics can now be rewritten in the affine form as 
\begin{equation}
	\m{\dot{x}}(t) = \m{a}(\m{x}(t),u(t)) = \m{g}(\m{x}(t)) + \m{h}(\m{x}(t))u(t),
\end{equation}
where
$\m{g}(\m{x}(t))$ and $\m{h}(\m{x}(t))$ are not functions of the control. The Hamiltonian  from Eq.~\eqref{eq:controlHamiltonian} is redefined as
\begin{equation}\label{eq:hamiltonian} 
  \C{H}(\m{x}(t),\g{\lambda}(t),u(t),t) = \m{f}(\m{x}(t),\g{\lambda}(t)) + \phi\T (\m{x(t)},\g{\lambda}(t))u(t),
\end{equation}
where $\m{f}(\m{x}(t),\g{\lambda}(t))$ and $\phi(\m{x}(t),\g{\lambda}(t))$ are the components of the Hamiltonian that are not a function of the control, and mixed state and control path constraints are not considered.
If the following holds along an optimal control
\begin{equation}\label{eq:Hu}
  \frac{\partial \C{H}}{\partial u} = \C{H}_{u} =  \phi(\m{x}(t),\g{\lambda}(t)) = 0,
\end{equation}
then the minimizing control in Eq.~\eqref{eq:PMP} is the interval $[u_{\min},u_{\max}]$; thus Eq.~\eqref{eq:PMP} only implies that an optimal control is feasible. Note that the control does not appear in Eq.~\eqref{eq:Hu} because the Hamiltonian is linear in the control. A singular arc is characterized as $\C{H}_u = 0$ and $\C{H}_{uu}$ is singular everywhere on the arc. When this occurs, the reduced Hessian matrix associated with the corresponding NLP that arises from the direct transcription method of Section~\ref{sect:LGR} is ill-conditioned such that the projected Hessian matrix is not positive definite.  This leads to poor conditioning in the control profiles which often presents itself in the form of oscillations, or chattering behavior, in the control solution.
 
The sign and value of $\phi(\m{x}(t),\g{\lambda}(t))$ (where $\phi$ is called the switching function) determines if the control is called a {\em bang-bang} control or a {\em singular} control. The weak form of PMP~\cite{Bryson1975} is used in the case of nonsmooth control and the minimization of the Hamiltonian leads to the following piecewise-continuous control, $u$, that is dependent on the switching function as follows
\begin{equation}\label{eq:optimal-control-singular}
  u^* = \operatorname*{arg\,min}_{u \in \, [u_{\min},u_{\max}]} \C{H}
    =\left\{
    \begin{array}{lcl}
      u_{\min} & , & \phi(\m{x}(t),\g{\lambda}(t))>0, \\
      u_{s} & , &\phi(\m{x}(t),\g{\lambda}(t))=0, \\
      u_{\max} & , & \phi(\m{x}(t),\g{\lambda}(t))<0,
    \end{array}
    \right.
\end{equation}
where the sign of the switching function, $\phi(\m{x}(t),\g{\lambda}(t))$, is determined by the state and the costate and $u_{s}$ lies in the closed interval $[u_{\min},u_{\max}]$. As the switching function $\phi(\m{x}(t),\g{\lambda}(t))$ changes sign, the control coincides with the sign changes by switching between its maximum and minimum values. Any time interval over which $\phi(\m{x}(t),\g{\lambda}(t))$ is zero is referred to as a singular arc and any control in the admissible control set will minimize the Hamiltonian. Furthermore, switching between nonsingular and singular arcs give rise to discontinuities on the state and control profiles, and the location of these transition points are referred to as {\em switch times}. These discontinuities defined by the switch times create numerical issues, whereas, the singular control suffers from non-uniqueness issues that occur when the control is free to lie between its upper and lower bounds and is not defined by the optimality conditions. 


The singular  control is obtained implicitly from the switching function. Specifically, $\phi$ is differentiated repeatedly until the control $u$ explicitly appears~\cite{Schattler2012}. Therefore, $u^*$ can be solved for by
\begin{equation}\label{eq:singularCondition}
	\frac{d^{(2r)}}{dt^{(2r)}} \phi= 0, \quad (r = 0,1,2,\ldots), 
\end{equation}
where $2r$ is the minimum number of differentiations of $\phi$ required to obtain the corresponding control $u_{s}$.
For $u$ to be optimal over a singular arc, the number of differentiations $2r$ must be even~\cite{Bryson1975,Schattler2012}.
Furthermore, the generalized Legendre-Clebsch condition~\cite{Bryson1975,Kelley1967,Kopp1965}
\begin{equation}\label{eq:generalized-Legendre-Clebsch-condition}
	(-1)^r \frac{\partial}{\partial u}\left[ \frac{d^{2r}}{dt^{(2r)}}\phi \right] \geq 0, \quad (r = 0,1,2,\ldots), 
\end{equation}
must hold over the duration of a singular arc.
While in some problems of interest it is possible to use Eq.~\eqref{eq:singularCondition} to determine a condition for the singular control, in many cases it is unable to produce the singular control (for example, if the order of the singular arc is infinite). Even in cases where the singular control could be determined from Eq.~\eqref{eq:singularCondition}, taking derivatives higher than second-order is not easy to implement numerically or analytically. As mentioned previously, the singular control might be a function of both the state and the costate. If this is the case then a direct collocation method could not be utilized to determine the optimal trajectory.  In this paper a method is developed that can be fully automated for detecting and accurately approximating the solution of bang-bang and singular optimal control problems.

\section{Method for Bang-Bang \& Singular Optimal Control Problems\label{sect:method}}

In this section the method for solving bang-bang and singular optimal control problems is developed. The method consists of two stages. The first stage of the method described in Section~\ref{sect:sructureDet} details the detection of the control structure and the decomposition of the optimal control problem into a multiple-domain optimal control problem dictated by the discontinuities identified that are represented as domain interface variables. These domain interface variables are then treated as additional decision variables in the nonlinear programming problem (NLP). The first stage is only implemented on the first mesh iteration. The second stage described in Section~\ref{sect:constraints} describes the new constraints that are added to the NLP depending on the structure detection's classification of a domain as being bang-bang or singular in order to constrain the modified optimal control problem correctly. The constraints and methods applied in each type of domain are provided in Section~\ref{sect:constraints}. The second stage initiates the iterative procedures in the proposed method described in Section~\ref{sect:algorithm}.

\subsection{Structure Detection and Decomposition\label{sect:sructureDet}}

Assume now that the optimal control problem formulated in Section~\ref{sect:probForm} under the assumptions of Section~\ref{sect:SOCP} has been transcribed into a NLP using multiple-domain LGR collocation developed in Section~\ref{sect:LGR} with $D=1$ (that is, a single domain is used).  The solution obtained from the NLP then leads to estimates of the state, the control, and the costate as given in Eqs.~\eqref{eq:NLP-Var-Mat} and~\eqref{eq:NLP-costates}, respectively.  Assume further that the mesh refinement accuracy tolerance is not satisfied.  As a result, mesh refinement is required which simultaneously enables the decomposition of the problem into domains that are either bang-bang, singular, or regular. This decomposition is obtained using structure detection as described now.  

Structure detection locates discontinuities identified on the initial mesh and then uses the locations of the discontinuities to determine the classification of the interval formed by two adjacent discontinuities. In this work, only control discontinuities are considered because only problems where the Hamiltonian is linear in control are analyzed. Their locations are estimated using jump function approximations~\cite{Miller2021} and then the intervals formed by each discontinuity is analyzed using the switching function. Structure detection begins by applying the method of Section~\ref{sect:jpDet} to identify and estimate the locations of any control discontinuities. After discontinuity locations have been estimated, the method of Section~\ref{sect:hamDet} takes the estimated discontinuity locations and determines the classification of the domain as bang-bang or singular. The structure detection and decomposition process only occurs once on the initial solution.  Detailed descriptions of the structure detection process are described next.

\subsubsection{Identification of Control Switch Times\label{sect:jpDet}}

Discontinuities in each component of the control are identified using jump function approximations of the control solution as shown in Ref.~\cite{Miller2021}.  In particular, the method given in Ref.~\cite{Miller2021} is employed here because it is effective for estimating locations of nonsmoothness in the optimal control.  A brief overview of the process given in Ref.~\cite{Miller2021} is provided here for completeness.  For further details related to jump function approximations for detecting nonsmoothness in an optimal control, see Ref.~\cite{Miller2021}.


First, a jump function is defined as follows. Let $f:\R\rightarrow\R$ be an arbitrary function defined on the interval $t\in[t_0,t_f]$. The jump function of $f(t)$, denoted $[f](t)$, is defined as $[f](t)=f (t^+)-f(t^-)$ where $f(t^+)$ and $f(t^-)$ are the right-hand and left-hand limits of $f(t)$, $t\in[t_0,t_f]$. The jump function is zero across intervals where $f(t)$ is continuous and takes on the value of the jump in $f(t)$ at those locations where $f(t)$ is discontinuous.
According to Refs.~\cite{Miller2021,Archibald2005} the jump function of a function $f(t)$ is approximated by
\begin{equation} \label{eq:Lmf}
	L_m f(t) = \frac{1}{q_m (t)} \sum_{t_j \in \C{S}_t} c_j(t) f(t_j) \approx [f](t),
\end{equation}
where $q_m (t)$ is defined by
\begin{equation} \label{eq:qm}
	q_m (t) = \sum_{t_j \in \C{S}_{t}^{+}} c_j(t),
\end{equation}
$c_j (t)$ is defined by
\begin{equation} \label{eq:cj}
	c_j (t) = \frac{m!}{\prod\limits_{\substack{i = 1 \\i \neq j}}^{m + 1} (t_j - t_i)},
\end{equation}
and $m$ specifies the order of the approximation.  Higher order approximations converge to the jump function faster outside the neighborhood of discontinuities but have oscillatory behavior in the vicinity of discontinuities.  The oscillations are reduced via the \textbf{minmod} function, defined here as
\begin{equation} \label{eq:MMLmf}
MM\left( L_{\C{M}} f(t) \right) = \left\{
  \begin{array}{ll}
    \min\limits_{m \in \C{M}} L_m f(t), &  L_m f(t) > 0~~\forall~m \in \C{M}, \\
    \max\limits_{m \in \C{M}} L_m f(t), & L_m f(t) < 0~~\forall~m \in \C{M}, \\
    0                           & \textrm{otherwise},\\
  \end{array}\right.
\end{equation}
where $\C{M} \subset \bb{N}^{+}$ is a finite set of choices of the approximation order $m$.  

Suppose an initial control solution is obtained. The control solution, $U(\tau_j^{(k)}), j=\{1,\ldots,N_k\};~k=\{1,\ldots,K\}$, is normalized to the interval $[0,1)$ by the transformation
\begin{equation}
	u(\tau_j^{(k)}) = \frac{U(\tau_j^{(k)})-u_{\min}}{1+u_{\max}-u_{\min}}~,
\end{equation}
where $u_{\min}$ and $u_{\max}$ are the minimum and maximum values of the control. Together, the normalized control solution and the corresponding collocation points, $\tau_j^{(k)}$ on $[-1,1)$ of the initial mesh, are applied to Eqs.~\eqref{eq:Lmf}--\eqref{eq:MMLmf} to produce a jump function approximation for the normalized control. The jump function approximation is then evaluated at the points $\tau_{j+\frac{1}{2}}^{(k)} = \frac{1}{2}(\tau_j^{(k)} + \tau_{j+1}^{(k)}),  j=\{1,\ldots,N_k\};~k=\{1,\ldots,K\}$. Let the evaluation of the jump function approximation at  $\tau_{j+\frac{1}{2}}^{(k)}$ be denoted by $MM(\tau_{j+\frac{1}{2}}^{(k)})$. The method detects a discontinuity at the location, $\tau_{j+\frac{1}{2}}^{(k)}$, if the following condition is satisfied:
\begin{equation}\label{eq:JFdisc}
	\abs{MM(\tau_{j+\frac{1}{2}}^{(k)})} \geq \eta.
\end{equation}
The identified discontinuities are referred to as $b_i, \, i = \{1,\ldots,n_d\}$ where $n_d$ is the total number of identified control discontinuities.
It is noted that $\eta\in[0,1)$ in Eq.~\eqref{eq:JFdisc} is a user-specified threshold that specifies the relative size of jumps that are detected (where the likelihood of jumps being detected in the control decreases as $\eta$ increases). The default value of $\eta$ is set to $0.1$ for the examples considered in this paper.

Bounds on the discontinuity locations are now defined.  Consider for some $j=\{1,\ldots,N_k\}$ and $k=\{1,\ldots,K\}$ that Eq.~\eqref{eq:JFdisc} is satisfied, indicating that a discontinuity is present somewhere on the mesh interval $\tau \in [\tau_j^{(k)},\tau_{j+1}^{(k)}]$. To account for the uncertainty incurred by using the numerical solution as a sample for the jump function approximation, a safety factor, $\mu \geq 1$ is introduced to extend the bounds estimated for the discontinuity. This safety factor provides a larger threshold to adequately capture the potential search space of the estimated switch times. Let $[b^-_i,b^+_i], \, i = \{1,\ldots,n_d\}$ be the lower and upper bounds on the locations of discontinuities in the control (that is, any discontinuity is bounded to lie on the interval $[b^-_i,b^+_i]$). The estimates of these bounds are defined as
\begin{equation}
  \left.
    \begin{array}{lcl}
      b_i^- &=& \tau_{j+\frac{1}{2}}^{(k)} - \mu\left(\tau_{j+\frac{1}{2}}^{(k)} - \tau_{j}^{(k)}\right), \\ \\
      b_i^+ &=& \tau_{j+\frac{1}{2}}^{(k)} + \mu\left(\tau_{j+1}^{(k)} - \tau_{j+\frac{1}{2}}^{(k)}\right),
    \end{array} \right\} \left. \begin{array}{ll} i & = \{1,\ldots,n_d\}, \\ j & =  \{1,\ldots,N_k\}, \\ k& =\{1,\ldots,K\}.\end{array} \right.
\end{equation}
Larger values of $\mu$ are more desirable as it is more likely that the discontinuity $b_i$ will lie in the interval $[b_i^-,b_i^+],~i=\{1,\ldots,n_d\}$. Furthermore, the default value of $\mu$ is set to $1.5$ for the examples considered in this paper.

\subsubsection{Identification of Bang-Bang and Singular Domains\label{sect:hamDet}}

Identification of the domains is employed to determine the structure of the control using the discontinuities identified from the jump function approximations. Specifically, the control solution is inspected to determine if any bang-bang or singular arcs exist. The first and second derivatives of the Hamiltonian with respect to control are computed and used to detect if the Hamiltonian is linear in the control by assessing if the second derivatives are zero. If the Hamiltonian is affine in control, the first derivatives are computed and represent the switching function of the system. If the Hamiltonian is not linear in the control the structure detection process is finished and smooth mesh refinement (see Section~\ref{sect:smooth}) can be performed.

Suppose the initial solution contains the following newly identified discontinuities $b_i$ and corresponding bounds $[b_i^-,b_i^+], \, i=\{1,\ldots,n_d\}$. The solution is divided into intervals starting with the initial time, the discontinuity locations, and ending with the final time, $\{[-1,b_1],[b_1,b_2],\ldots,[b_i,b_{i+1}],$ $[b_{i+1},+1]\}, \, i = \{1,\ldots,n_d\}$. Next, the Hamiltonian in Eq.~\eqref{eq:hamiltonian} is computed using the initial solution
\begin{equation}\label{eq:hamiltonianNLP} 
  \C{H} = \m{f}(\m{X}(\tau),\g{\lambda}(\tau)) + \phi\T (\m{X(\tau)},\g{\lambda}(\tau))U(\tau).
\end{equation}
Note that the costates are also obtained when solving the NLP that results from multiple-domain LGR collocation (see Section~\ref{sect:costates}). The first and second derivatives with respect to the control are computed using the already computed derivatives required by the NLP solver
\begin{eqnarray}
	\frac{\partial \C{H}}{\partial U} = \phi(\m{X}(\tau),\g{\lambda}(\tau),t(\tau,t_0,t_f)),  \label{eq:firstHu}  \\
	\frac{\partial^2 \C{H}}{\partial U^2} = \Phi(\m{X}(\tau),\g{\lambda}(\tau),t(\tau,t_0,t_f)).  \label{eq:secondHu}
\end{eqnarray}
First, the values of Eq.~\eqref{eq:secondHu} must be zero. If this condition is satisfied then Eq.~\eqref{eq:firstHu} is analyzed as follows.
A bang-bang interval in the control structure will occur when the switching function $\phi$ changes sign. The sign of the  switching function $\phi$ is checked in each interval $[b_i,b_{i+1}], \, i = \{1,\ldots,n_d\}$. If the sign in the interval is positive, the control is constrained to its minimum value. If the sign in the interval is negative, the control is constrained to its maximum value. Additional details on these constraints are discussed in Section~\ref{sect:bang}.

A singular interval in the control structure will occur when the switching function $\phi$ is zero at every point in the interval $[b_i,b_{i+1}]$. Assessing if the switching function $\phi$ is zero over the current interval is not a trivial task. Due to the NLP being ill-conditioned when a singular arc is present, the estimated solution over a singular interval will suffer from larger numerical error. The user-defined zero threshold becomes critical in detecting the presence of a singular arc because the switching function will never be exactly zero. This threshold is heavily influenced by the coarseness of the initial mesh and the accuracy of the detected discontinuity locations. Once a singular arc has been detected, a regularization method is employed as described in Section~\ref{sect:singular} and~\ref{sect:iterativeSingular}.  If a scenario occurs where the entire control is singular on $[t_0,t_f]$, then no discontinuities will be detected. In this situation, the identification procedure is applied over the entire time domain so that the singular arc can be identified.

\subsubsection{Structure Decomposition\label{sect:structureDecomp}}

Assuming the methods of Section~\ref{sect:jpDet} and~\ref{sect:hamDet} have identified discontinuities and intervals that are bang-bang or singular, the initial mesh is now decomposed into the multiple-domain structure. Once acquired, the detected structure of the nonsmooth control is used to introduce the appropriate number of domain interface variables, $t_s^{[d]}, \, d = \{1,\ldots,D-1\}$, to be solved for on subsequent mesh iterations, where the initial guess for each variable is the estimated discontinuity location $b_i,\,i=\{1,\ldots,n_d\}$ that was found using the method in Section~\ref{sect:jpDet}. The domain interface variables are included in the NLP by adding them as additional decision variables that define the new domains, $\C{P}_d=[t_s^{[d-1]},t_s^{[d]}], \, d = \{1,\ldots,D\} $. Specifically, the domain interface variables are employed by dividing the time horizon $t \in [t_0,t_f]$ of the original optimal control problem into $D$ domains as described in Section~\ref{sect:LGR}.  

Next, bounds on the domain interface variables are enforced to prevent the collapse or overlap of domains. The bounds provide an additional constraint on the domain interface variables. The upper and lower bounds on each domain interface variable are determined by taking the discontinuity bounds found in Section~\ref{sect:jpDet} and transforming them to the time interval $t \in [t_0,t_f]$ using the transformation in Eq.~\eqref{eq:affine-transformation}. Thus, the bounds $[b_i^-,b_i^+], \, i = \{1,\ldots,n_d\}$ are transformed to $[t_l^{[d-1]},t_u^{[d]} ], \, d = \{1,\ldots,D-1\}$.  

\begin{figure}
\begin{center}
	\includegraphics[scale=0.8]{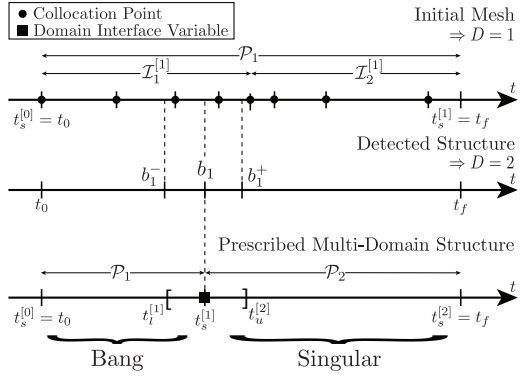}
	\caption{Schematic of process for decomposing the nonsmooth optimal control problem into $D$ domains where the $D-1$ domain interface variables are included as optimization variables to determine the optimal switch times in the control.\label{fig:structureDetection}}
\end{center}
\end{figure}

This approach to structure decomposition partitions the entire problem domain into multiple domains of the form described in Section~\ref{sect:LGR} such that the switch times are represented by the strategically placed domain interface variables $t_s^{[d]},\, d \in \{1,\ldots,D-1\}$.  A schematic for the process of decomposing the nonsmooth control structure into a multiple-domain formulation with domain interface variables is shown in Fig.~\ref{fig:structureDetection}. Additionally, the form of the control in each domain is classified as either bang-bang, singular, or regular.  In the next section, the constraints and refinement strategies required by each type of domain are discussed.

\subsection{Domain Constraints and Refinement}\label{sect:constraints}

Now that structure detection and decomposition has taken place by the methods of Section~\ref{sect:sructureDet}, additional constraints are required to properly constrain the multiple-domain optimal control problem. Recall that there are three types of domain classifications: bang-bang, singular, and regular. Each domain type requires its own set of constraints and refinement methods that are detailed in the following sections.

\subsubsection{Bang-Bang Domain Constraints\label{sect:bang}}

Suppose the problem has been partitioned into $D$ domains based on the results of structure detection, and it has been determined that $\C{B}$ domains are bang-bang by the method of Section~\ref{sect:hamDet}, where $\C{B}\leq D$.  The value of the switching function $\phi$ in the current domain of interest $\C{P}_d$ is used to determine the value of the control over that domain. 
Recall, that the switching function assigns the value of the control according to
\begin{equation}
	\begin{array}{lclcl}
		u^{[d]}(\tau) &=& u_{\min} &,& \phi^{[d]}(\m{X},\g{\lambda},t(\tau,t_0,t_f)) > 0, \\
		u^{[d]}(\tau) &=& u_{\max} &,& \phi^{[d]}(\m{X},\g{\lambda},t(\tau,t_0,t_f)) < 0. 
	\end{array}
\end{equation} 
The control is then constrained to its corresponding maximum or minimum value in the resulting multiple-domain optimal control problem. 
The bang-bang control is now appropriately constrained over its domain and the corresponding domain interface variables can be optimized to the optimal switch time locations.

\subsubsection{Regularization of a Singular Domain\label{sect:singular}}

Assume now that the entire domain $[t_0,t_f]$ of the optimal control problem has been partitioned into $D$ domains using the structure detection method as described in Section \ref{sect:sructureDet}.  Assume further that, using the procedure given in Section~\ref{sect:hamDet}, $\C{S}$ of these $D$ domains are classified as singular (where $\C{S}\leq D$ and $\C{S}+\C{B} \leq D$) such that $\{s_1,\ldots,s_\C{S}\}\subseteq \{1,\ldots,D\}$ are the indices corresponding to the singular domains.  The singular domains are then defined, respectively, on the intervals $[t_s^{[s_d-1]},t_s^{[s_d]}] \subseteq [t_0,t_f]$, $d=\{1,\ldots,\C{S}\} $.

In any domain that is classified as singular, the following iterative regularization method is employed.  First, the objective functional in the singular domain is augmented with the regularization term
\begin{equation}\label{eq:penaltyTerm}
  \delta_{s_d} \, = \, \frac{\epsilon}{2} \int_{t_s^{[s_d-1]}}^{t_s^{[s_d]}}  \left(u(t)-\alpha_{p}(t)\right)^2\, dt, \quad d = \{1,\ldots,S\},\, p = \{1,2,\ldots\}
\end{equation}
where $u(t)$ is the optimal control to be determined when solving the problem, $\alpha_p(t) \in \R$ is a known function that changes with each iteration of the regularization method (see Section \ref{sect:iterativeSingular}), and 
\begin{equation}\label{eq:penaltyTerm-integrand}
  \frac{\epsilon}{2} \left(u(t)-\alpha_{p}(t)\right)^2
\end{equation}
is the integrand of Eq.~\eqref{eq:penaltyTerm}.  Furthermore, $\epsilon$ is a user-defined weighting parameter that is chosen based on the particular problem under consideration.  Augmenting the Hamiltonian with the term in Eq.~\eqref{eq:penaltyTerm-integrand} results in a Hamiltonian that is quadratic in the control.  Consequently, the optimal control problem, that would be singular without the inclusion of the term in Eq.~\eqref{eq:penaltyTerm-integrand}, becomes regular (nonsingular).  While in principle $\epsilon$ can be any positive value, it must be sufficiently large to eliminate the indeterminacy of determining the optimal control on the singular arc, but it must be sufficiently small so that the resulting optimal control is in close proximity to the true singular optimal control.

\subsubsection{Iterative Procedure for Determining Singular Control\label{sect:iterativeSingular}}

Next, it is important to understand the source of the function $\alpha_p(t)$.  Singular domain refinement is employed in an iterative fashion where $p$ is the iteration number of the singular domain refinement procedure.  For $p=1$, $\alpha_p(t)$ is set to zero.  Then, for $p>1$, $\alpha_p(t)$ is obtained using a continuous piecewise cubic approximation of the control obtained from the solution of the NLP on iteration $p-1$ and the approximation for $\alpha_p(t)$ is obtained by interpolation using a piecewise cubic polynomial with the properties defined in Ref.~\cite{Fritsch1980}.  

\begin{figure}
\begin{center}
	\includegraphics[scale=0.6]{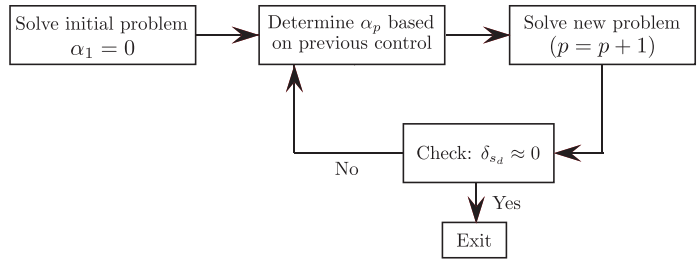}
	\caption{Description of the regularization method during mesh refinement iteration, $M$.\label{fig:proximalGraphic}}
\end{center}
\end{figure}

The augmented multiple-domain optimal control problem that is solved by the regularization method is then stated as follows.  First, for any iteration $p\geq1$ the objective functional to be minimized includes the terms $\delta_{s_d},~\{d=1,\ldots,S\},$ and is defined as
\begin{equation}\label{eq:aug-bolza-cost}
  \C{J}_a   = \C{J} +  \sum_{d=1}^{\C{S}} \delta_{s_d}.
\end{equation}
Furthermore, the constraints include the dynamic constraints, the boundary conditions, and the path constraints given, respectively, in Eqs.~\eqref{eq:single-dyn-t}--\eqref{eq:single-bc-t}.  The process for  updating the iteration of the regularization method is as follows.  First, after solving the NLP arising from LGR collocation with $D=1$, the structure detection and decomposition method of Section~\ref{sect:sructureDet} is employed for one mesh iteration.  On the second mesh iteration (that is, $M=2$) the value of $p$ is set to unity (that is, $p\equiv 1$) and the regularization term of Eq.~\eqref{eq:penaltyTerm} is augmented to the objective functional in any domain classified as singular with $\alpha_p(t)\equiv 0$ and the resulting NLP is then solved.  Then, from the third mesh iteration onwards (that is, $M>2$), the value of $p$ is incremented as $p\rightarrow p+1$ and the value of $\alpha_{p}(t)$ is determined based on the aforementioned continuous piecewise cubic approximation from the discrete control solution obtained on iteration $p$.  The process of refining the mesh and updating the function $\alpha_p(t)$ is then repeated until the regularization term $\delta$ lies below a user-specified tolerance $\sigma$. The regularization procedure is independent of regular mesh refinement meaning that refinement of the mesh always occurs on each iteration and the updates to $\alpha_p(t)$ occur independently from the mesh refinement actions. Figure~\ref{fig:proximalGraphic} provides a schematic of the regularization method that is enforced over each mesh refinement iteration in a singular domain.


\subsubsection{Regular Domain Refinement\label{sect:smooth}}

Suppose that a particular domain $\C{P}_d$ has been determined to be regular, that is, the domain $\C{P}_d$ is categorized neither as bang-bang nor singular.   Although the domain $\C{P}_d$ is categorized as regular, this domain may still require mesh refinement.  In particular, mesh refinement of the domain $\C{P}_d$ will be required if the maximum relative error on the domain $\C{P}_d$ exceeds the relative error tolerance.  
Several methods have been previously developed for smooth mesh refinement~\cite{Patterson2015,Liu2015,Darby2010,Liu2018}.  In this research, domains where the solution is smooth are refined using the method developed in Ref.~\cite{Patterson2015}.  
The mesh refinement method of Ref~\cite{Patterson2015} occurs on every iteration until the user-specified mesh error tolerance is met and is independent of the status of the regularization procedure in Section~\ref{sect:singular}.
It is noted that Ref.~\cite{Patterson2015} also includes a method for the computation of a relative error estimate.  The reader is referred to Ref.~\cite{Patterson2015} for a more detailed explanation of how both smooth mesh refinement operates and an estimate of the relative error on a mesh.  


\subsection{Procedure for Solving Bang-Bang and Singular Optimal Control Problems \label{sect:algorithm}}
An overview of the proposed method for bang-bang and singular optimal control problems is shown below. The mesh refinement iteration is denoted by $M$ and is incremented by one with each loop of the method. The regularization refinement iteration is denoted by $p$ and is incremented by one on each mesh iteration as required by the regularization method. This method terminates when two requirements are met. First, the regularization term \eqref{eq:penaltyTerm} must be within a user specified tolerance $\sigma$ from zero, or if Eq.~\eqref{eq:penaltyTerm} remains identical in value over three consecutive iterations. Second, the mesh error tolerance, $e$, must be satisfied on each mesh interval or if $M$ reaches a prescribed limit, $M_{\max}$. The method is executed as follows:

\begin{shadedframe}
\vspace{-10pt}
\begin{center}
 \shadowbox{\bf Method for Solving Bang-Bang \& Singular Optimal Control Problems} 
\end{center}
\begin{enumerate}[{\bf Step 1:}]
\item Set $M=0$ and specify initial mesh.  All mesh intervals form a single domain.
\item Solve NLP of Section~\ref{sect:LGR} on mesh $M$. \label{step:solve}
\item Compute relative error $e$ on current mesh $M$. \label{step:error}
\item If $M=1$, employ structure detection and decomposition in Section~\ref{sect:sructureDet}.
	\begin{enumerate}[{\bf (a):}]
	\item Determine the number of switch times, $n_d$, using the methods of Section~\ref{sect:jpDet}. \label{step:jpDet}
	\item Classify the intervals as bang-bang or singular using the method of~\ref{sect:hamDet}. \label{step:hamDet}
	\item Assign domain interface variables by method of Section~\ref{sect:structureDecomp}.
	\item Partition time horizon into domains by method of Section~\ref{sect:structureDecomp}.
	\item Perform domain refinement according to Section~\ref{sect:constraints}.
		\begin{enumerate}[{\bf (i):}]
		\item Enforce control constraints in each bang-bang domain as in Section~\ref{sect:bang}.
		\item Employ regularization method in singular domains according to Section~\ref{sect:singular}.
		\item Apply mesh refinement in regular domains as in Section~\ref{sect:smooth}.
		\end{enumerate}
	\end{enumerate}
\item If $n_d=0$ or $M>1$, apply smooth mesh refinement and proceed to {\bf Step \ref{step:end}}. 	
\item If ($M>1$, $\delta \leq \sigma$, and $e_{\max} \leq e$) or ($M > M_{\max}$), then quit. Otherwise: \label{step:checkProx}
	\begin{enumerate}[{\bf (a):}]
		\item Apply iterative procedure of Section~\ref{sect:iterativeSingular}.
	     \item Increment $M \rightarrow M+1,\, p \rightarrow p+1$ and return to {\bf Step \ref{step:solve}}
	\end{enumerate}
	\item Increment $M \rightarrow M+1$ and return to {\bf Step \ref{step:solve}}.\label{step:end}
\end{enumerate}
\end{shadedframe}


\section{Examples \label{sect:examples}}

In this section, the aforementioned BBSOC method described in Section~\ref{sect:method} is studied on three nontrivial bang-bang and singular optimal control problems. Each of the three problems demonstrates the methods abilities to solve different types of nonsmooth optimal control problems including a purely bang-bang control, a bang-singular-bang control with no analytic solution, and a bang-singular control where the analytic solution exists. A fourth example with a smooth optimal solution is also solved using the method of Section~\ref{sect:method} and compared against a previously developed mesh refinement method.  In addition, the following values of design parameters in the method are chosen for all examples: $\eta=0.1$ and $\mu=1.5$ while the value of the regularization parameter $\epsilon$ is different for each example and ranges between unity and $10^{-8}$ (depending upon the problem being solved).

All results obtained using the BBSOC method are compared with either an analytic solution (if an analytic solutions exists), a numerical solution obtained without enforcing the known structure of the optimal control, or a highly accurate numerical solution obtained by enforcing the known structure of the optimal control.  Any numerical solution obtained by enforcing the known structure of the optimal control is referred to as a {\em baseline} solution.  All numerical solutions, other than those obtained using the BBSOC method, are obtained using the MATLAB${}^{\textregistered}$ optimal control software $\mathbb{GPOPS-II}$ as described in Ref.~\cite{Patterson2014}, and $\mathbb{GPOPS-II}$ is referred to from this point forth as the $hp$--LGR method.  For any results obtained where the known structure of the optimal control is not enforced, the problem is formulated as a single-phase optimal control problem for use with the $hp$--LGR method.  For any results obtained where the known structure of the optimal control is enforced, the problem is formulated as a multiple-phase optimal control problem for use with $hp$--LGR method, and the control in each phase is either free (if the control is regular), is set to either its known lower or upper limit (if the control is bang-bang), or is determined by enforcing the singular arc optimality conditions (if the control is singular), and the switch times (which are the endpoints of each phase) are determined as part of the optimization.  Finally, for all numerical results obtained, the accuracy of the solution is improved using the mesh refinement method of Ref.~\cite{Patterson2015} where an error analysis is performed using the error estimate described in Ref.~\cite{Patterson2015}.   For completeness, Table~\ref{tab:different-methods} provides a table with the nomenclature that identifies the various methods being compared in this section.

\begin{table}[h]
  \centering
  \caption{Nomenclature of the various methods being compared in Section \ref{sect:examples}.\label{tab:different-methods}}
  \begin{tabular}{|c|c|}\hline
    {\bf Method Name} & {\bf Meaning} \\ \hline \hline
    BBSOC & Method Developed in This Paper \\ \hline
    $hp$-LGR & One-Phase Implementation of $\mathbb{GPOPS-II}$ Without Enforcement of Control Structure \\ \hline
    Baseline & Multiple-Phase Implementation of $\mathbb{GPOPS-II}$ With Enforcement of Control Structure \\ \hline
  \end{tabular}
\end{table}

All results are obtained using MATLAB${}^{\textregistered}$ and the nonlinear program developed in Section~\ref{sect:LGR} is solved using IPOPT~\cite{Biegler2008} in full-Newton mode. The NLP solver tolerance is set to $10^{-8}$ and the mesh refinement tolerance for smooth mesh refinement~\cite{Patterson2015} is $10^{-6}$.  All first and second derivatives are supplied to IPOPT using the automatic differentiation software ADiGator~\cite{Weinstein2017}.  In each example, the initial mesh consists of ten uniformly spaced mesh intervals and four collocation points per mesh interval, and the initial guess for all examples is a straight line for variables with boundary conditions at both endpoints and is a constant for variables with boundary conditions at only one endpoint.  Finally, all computations were performed on a 2.9 GHz 6-Core Intel Core i9 MacBook Pro running Mac OS Big Sur Version 11.3 with 32 GB 2400 MHz DDR4 of RAM, using \textsf{MATLAB} version R2019b (build 9.7.0.1190202) and all computation (CPU) times are in reference to this aforementioned machine.

\subsection*{Example 1: Robot Arm Problem}

Consider the following problem where the goal is to reorient a robotic arm in minimum time~\cite{Dolan2004}:
\begin{equation}\label{eq:robotArm}
  \begin{array}{lcl}
    &\textrm{minimize}& \,\, \C{J} = t_f, \vspace{0.5cm}  \\
    &\textrm{subject to}&
    \left\{
    \begin{array}{lclclclclcl}
      \dot{y_1}(t) &=& y_2(t)     &,&  y_1(0) &=& 9/2 &,& y_1(t_f) &=& 9/2, \\
      \dot{y_2}(t) &=& u_1(t)/L  &,&  y_2(0) &=& 0    &,& y_2(t_f) &=& 0, \\
      \dot{y_3}(t) &=& y_4(t)     &,&  y_3(0) &=& 0    &,& y_3(t_f) &=& 2\pi/3, \\
      \dot{y_4}(t) &=& u_2(t)/I_{\theta} &,&   y_4(0) &=& 0      &,& y_4(t_f) &=& 0, \\
      \dot{y_5}(t) &=& y_6(t)                 &,&   y_5(0) &=& \pi/4 &,& y_5(t_f) &=& \pi/4, \\
      \dot{y_6}(t) &=& u_3(t)/I{\phi}      &,&    y_6(0) &=& 0     &,& y_6(t_f) &=& 0, \\
      -1 &\leq& u_i(t) \leq 1, \quad (i = 1,2,3) &,&
    \end{array}
    \right.
  \end{array}
\end{equation} 
where $t_f$ is free, $I_{\phi} = \frac{1}{3}((L-y_1(t))^3+y_1^3(t))$, $I_{\theta} = I_{\phi}\textrm{sin}^2(y_5(t))$, and $L = 5$. The robot arm problem has a bang-bang structure for all three components of the control and contains a total of five switch points in the control. A baseline solution is obtained for comparison and the results are provided in Table~\ref{tab:robotArm}.

\begin{figure}
\centering
\vspace*{0.25cm}
\subfloat[\label{fig:robotArm-U1}]{\includegraphics[scale=0.37]{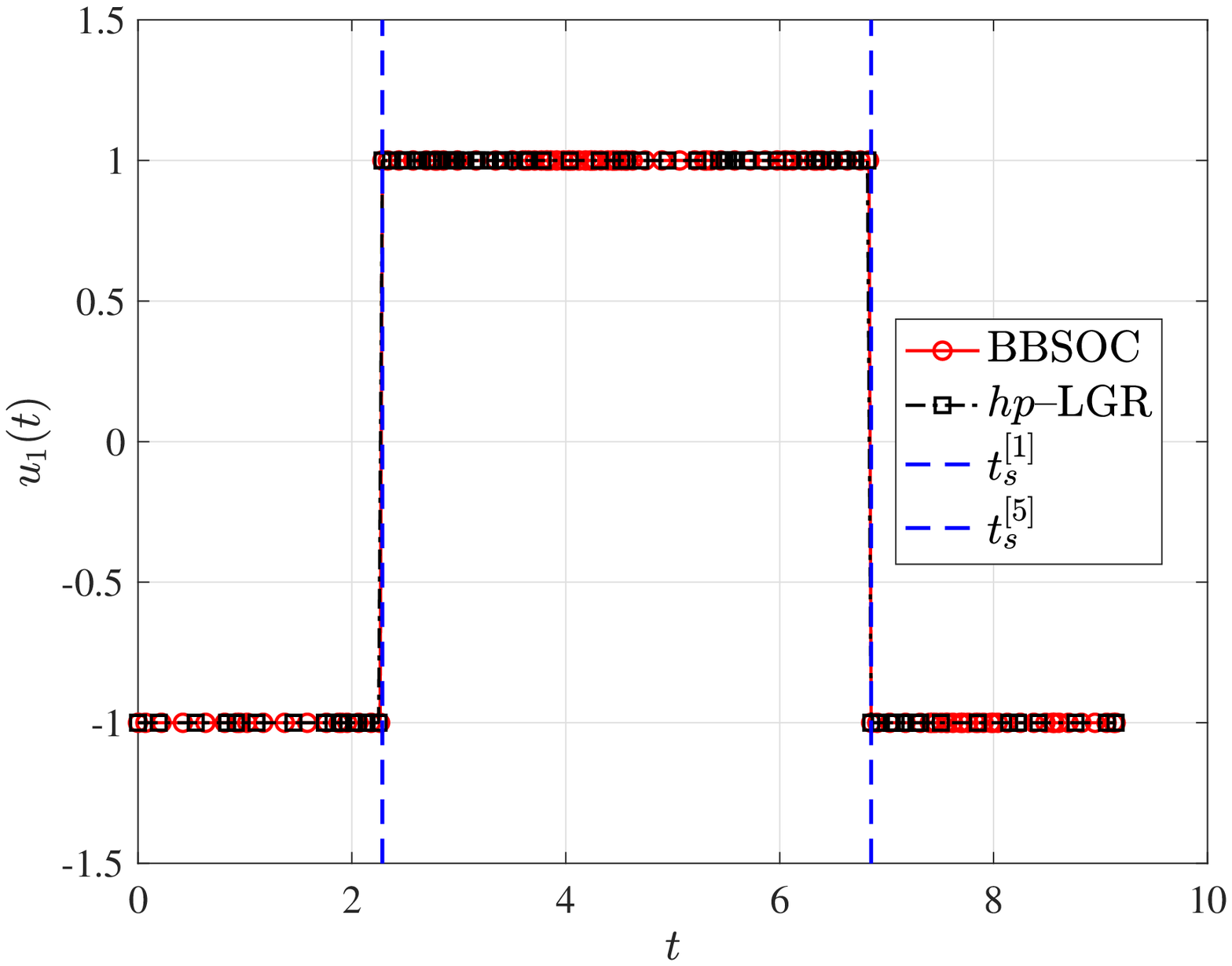}}~~\subfloat[\label{fig:robotArm-U2}]{\includegraphics[scale=0.37]{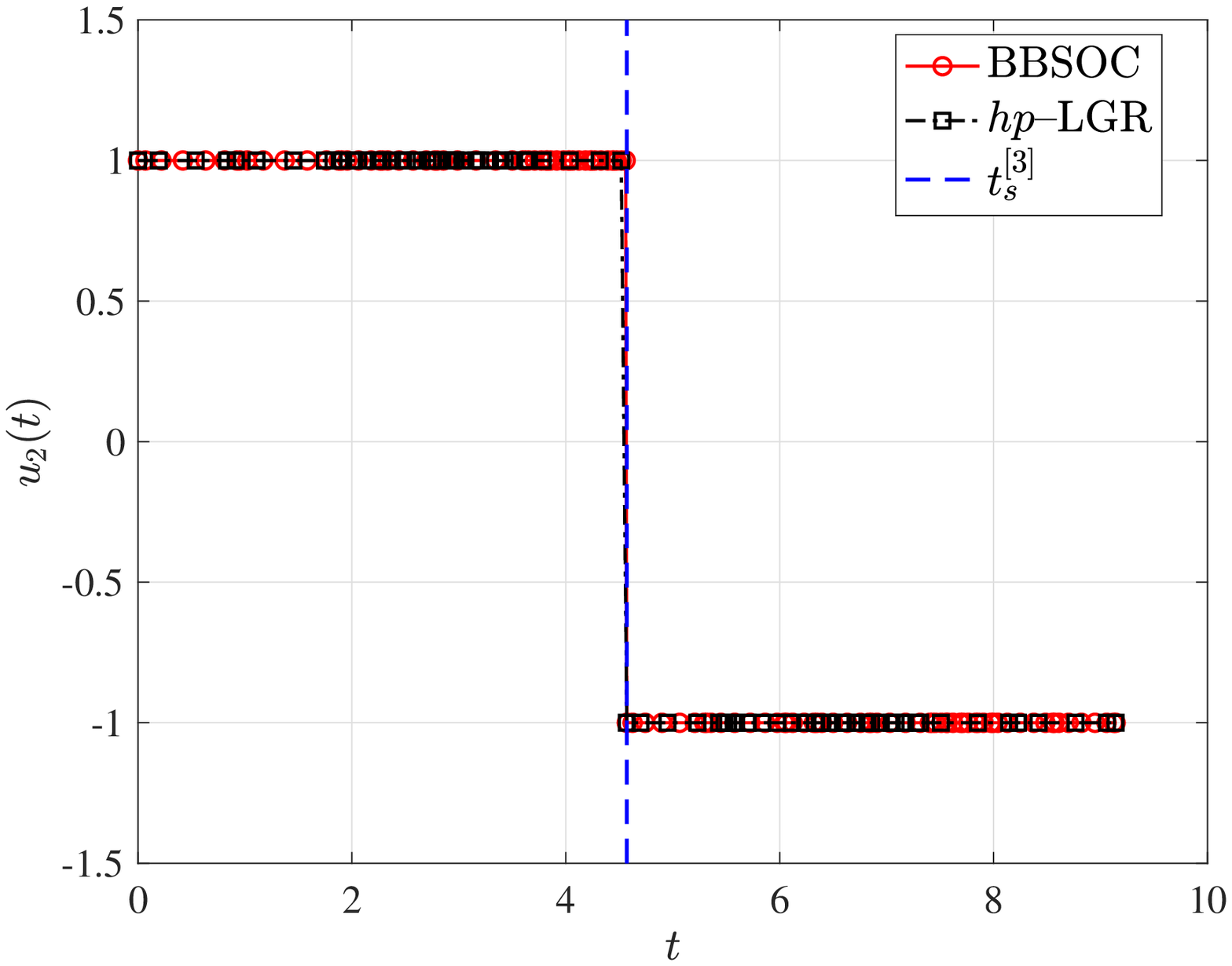}} \\
\subfloat[\label{fig:robotArm-U3}]{\includegraphics[scale=0.37]{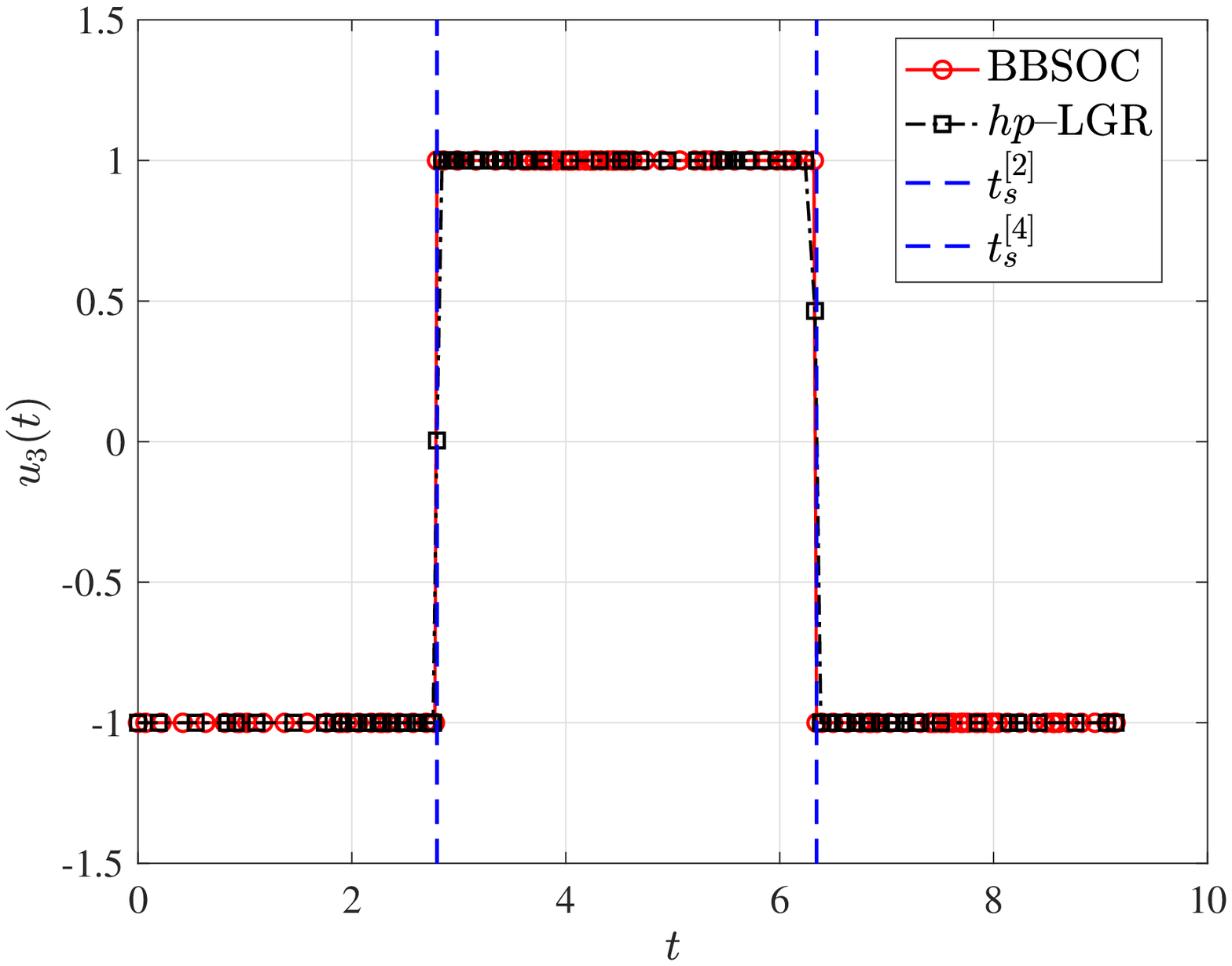}}
\caption{Control component solutions for Example~$1$ illustrate the discrepancies between the solutions obtained by the BBSOC and $hp$--LGR methods.} 
\label{fig:robotArm}
\end{figure}

The control solution obtained from solving this problem using the BBSOC and the $hp$--LGR method are shown in Fig.~\ref{fig:robotArm}.  Observing the controls in Figs.~\ref{fig:robotArm-U1}--\ref{fig:robotArm-U3}, it is seen that the five switch times are identified and the controls are constrained to their corresponding boundaries after just one iteration of the BBSOC method. In contrast, the $hp$--LGR method does not correctly identify all five discontinuities and as a result obtains incorrect approximations of the controls.  In particular, the $hp$--LGR method attempts to place more collocation points in the neighborhood of the switch times (thus increasing the size of the NLP) while still does not correctly identifying the exact locations of the switches in the control. Table~\ref{tab:robotArm} shows the numerical values of the switch times obtained from the BBSOC method, $hp$--LGR method, and the baseline solutions.  It is seen that the BBSOC method produces highly accurate approximations of both the switch times and objective and is in excellent agreement with the baseline solution.  On the other hand, the $hp$--LGR method (which is not designed specifically for solving problems with nonsmooth solutions) attains a less accurate solution when compared with the BBSOC method.  In addition to accuracy, Table~\ref{tab:robotArm} provides a comparison of the CPU time required for each approach.  In particular, it is seen from Table~\ref{tab:robotArm} that the BBSOC method converges to the optimal solution more efficiently compared with either the $hp$--LGR method or the baseline solution.  

\begin{table}
\caption{\label{tab:robotArm} Comparison of computational results for Example~$1$.}
\centering
\begin{tabular}{ c " c c c c c c c } 
\hline \hline
                     &  $t_{s}^{[1]}$   & $t_{s}^{[2]}$ &   $t_{s}^{[3]}$  & $t_{s}^{[4]}$  & $t_{s}^{[5]}$  & $\C{J}^*$ & CPU~$[\textrm{s}]$ \\ \thickhline 
BBSOC     & $2.285228$  & $2.796043$  & $4.570456$  & $6.344869$  &  $6.855684$ & $9.140912$ & $0.63$ \\ 
$hp$--LGR   & $2.249813$   & $2.795031$  & $4.518293$  & $6.330691$  & $6.820305$  & $9.140984$ & $3.25$  \\
Baseline    & $2.285228$ & $2.796043$ & $4.570456$ & $6.344869$ & $6.855684$ & $9.140912$ & $1.59$ \\
\hline \hline
\end{tabular}
\end{table}

\subsection*{Example 2: Goddard Rocket Problem}

Consider the following optimal control problem \cite{Bryson1975}:
\begin{equation}\label{eq:goddardRocket}
	\begin{array}{lcl}
    &\textrm{minimize}& \,\, \displaystyle \C{J} =  -h(t_f), \vspace{0.5cm} \\
    &\textrm{subject to}& 
    \left\{
    \begin{array}{lclclclclcl}
		\dot{h}(t) & = & v(t)                            &,&  h(0) & = & 0 & , & h(t_f) & = & \textrm{Free}, \\
		\dot{v}(t) & = & \frac{T(t)-D}{m(t)} - g    &,&  v(0) & = & 0 & , & v(t_f) & = & \textrm{Free}, \\
		\dot{m}(t) & = & -\frac{T(t)}{c}            &,&  m(0) & = & 3 & , & m(t_f) & = & 1, \\
		0 &\leq& T(t) \,\,\, \leq T_{\max} &,&
	\end{array}
    \right.
    \end{array}
\end{equation} 
where $h$ is the altitude, $v$ is the velocity, $m$ is the mass, $T$ is the thrust (and is the control), $D = D_0 v^2(t) \exp(-h(t)/H)$, and the final time is free. Further details on the model and the parameters can be found in Ref.~\cite{Bryson1975}.  In particular, the control $T$ has a bang-singular-bang structure and differentiating $\C{H}_u$ twice with respect to time leads to the following singular control law \cite{Betts2010}: 
\begin{equation}\label{eq:GR-singularControl}
  T_{\textrm{sing}}(t) = D + m(t)g + \left[ \frac{c^2(1+\frac{v(t)}{c})}{Hg} -1 - \frac{2c}{v(t)} \right ] \left[ \frac{m(t)g}{1 + \frac{4c}{v(t)} + \frac{2c^2}{v^2(t)}} \right ].
\end{equation}
The singular arc condition in Eq.~\eqref{eq:GR-singularControl}is used to obtain a baseline solution and the results are included in Table~\ref{tab:goddardRocket} for comparison.

\begin{figure}
\centering
\vspace*{0.25cm}
\subfloat[\label{fig:goddardRocket-U}]{\includegraphics[scale=0.4]{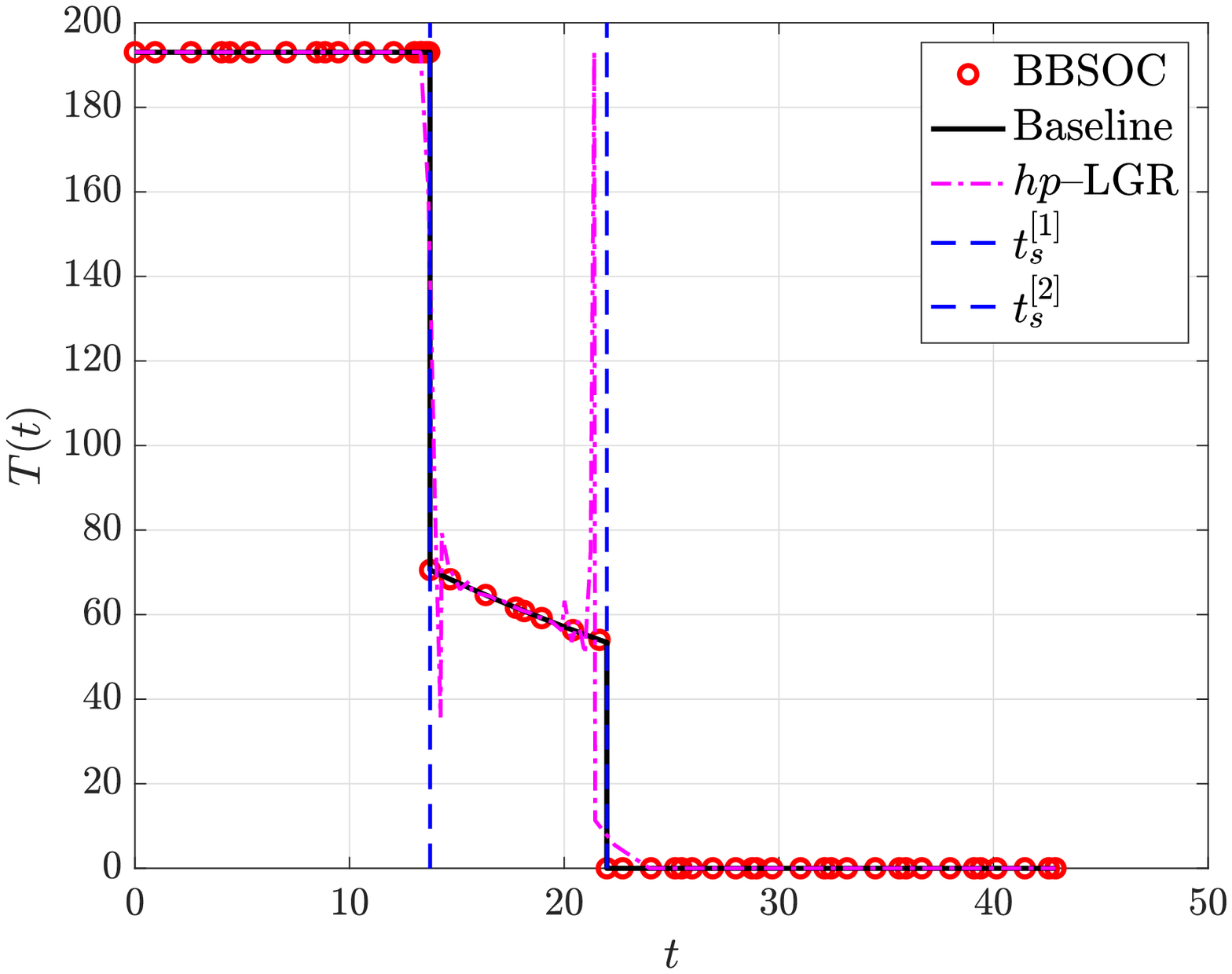}}~~\subfloat[\label{fig:goddardRocket-Hist}]{\includegraphics[scale=0.4]{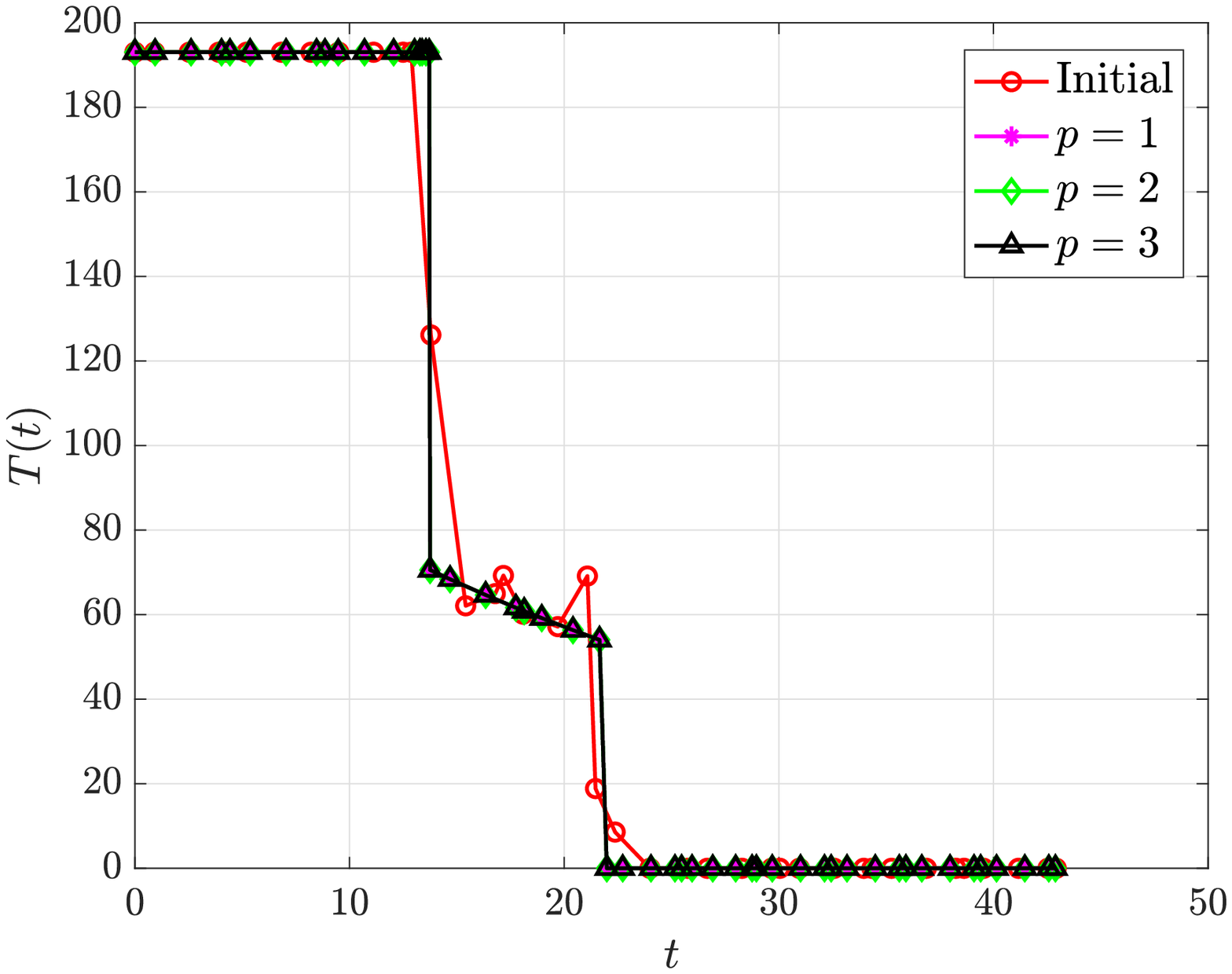}} \\ 
\caption{Control solution for Example~$2$ in Fig.~\ref{fig:goddardRocket-U} illustrates the numerical discrepancies between the solutions obtained by the BBSOC method, $hp$--LGR method, and the baseline method. The control history of the regularization procedure shows the process of obtaining the singular control in Fig.~\ref{fig:goddardRocket-Hist}.} 
\label{fig:goddardRocket}
\end{figure}

The problem in Eq.~\eqref{eq:goddardRocket} is solved using the BBSOC and $hp$--LGR methods. Both of these control solutions are compared in Fig.~\ref{fig:goddardRocket-U} and parameters related to the regularization method are provided in Table~\ref{tab:goddardRocket}. Figure~\ref{fig:goddardRocket-U} demonstrates that the regularization method has converged to the correct singular control while also correctly identifying the switch times.   In contrast, the $hp$--LGR method obtains a solution that exhibits oscillations at the switch times defining the singular arc. These oscillations are also seen in the initial iterations of the regularization method, but are removed by the third iteration as seen in Fig.~\ref{fig:goddardRocket-Hist}.

Next, Table~\ref{tab:goddardRocket} compares the solution obtained using the BBSOC method, $hp$--LGR method, and the baseline solution. These results further validate the BBSOC methods ability to identify the optimal switch times, final time, and cost. The BBSOC results are also in close agreement with the baseline solution while the $hp$--LGR results are not. The computation times for all three methods are also compared in Table~\ref{tab:goddardRocket} and all three methods have very similar CPU times with the BBSOC method taking the longest; however, it produces a far more accurate solution and achieves the same accuracy as the baseline solution without any a priori knowledge of the singular problem.

\begin{figure}
\centering
\vspace*{0.25cm}
\subfloat[\label{fig:goddardRocket-Surf}]{\includegraphics[scale=0.4]{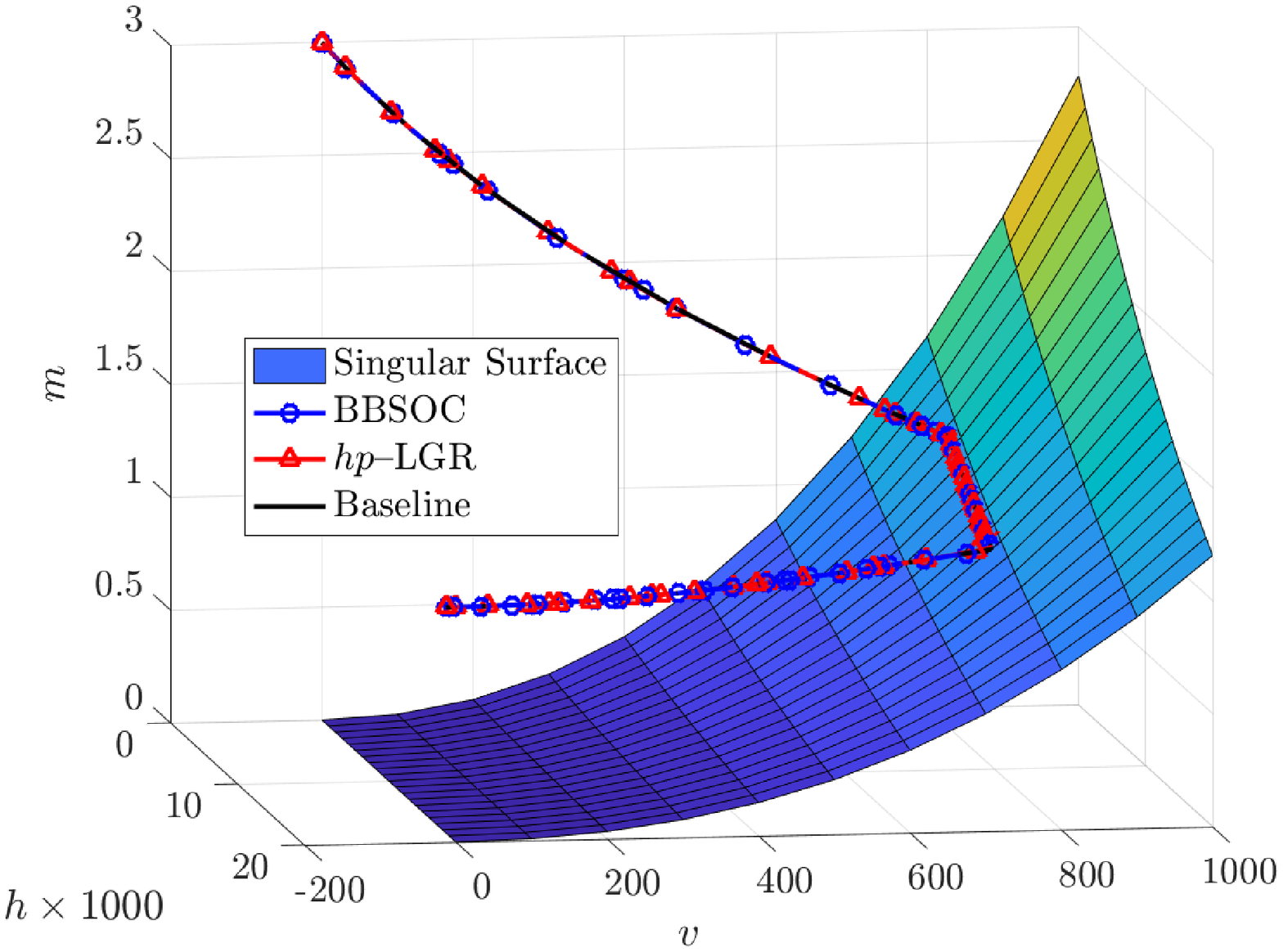}}~~\subfloat[\label{fig:goddardRocket-SurfZoom}]{\includegraphics[scale=0.4]{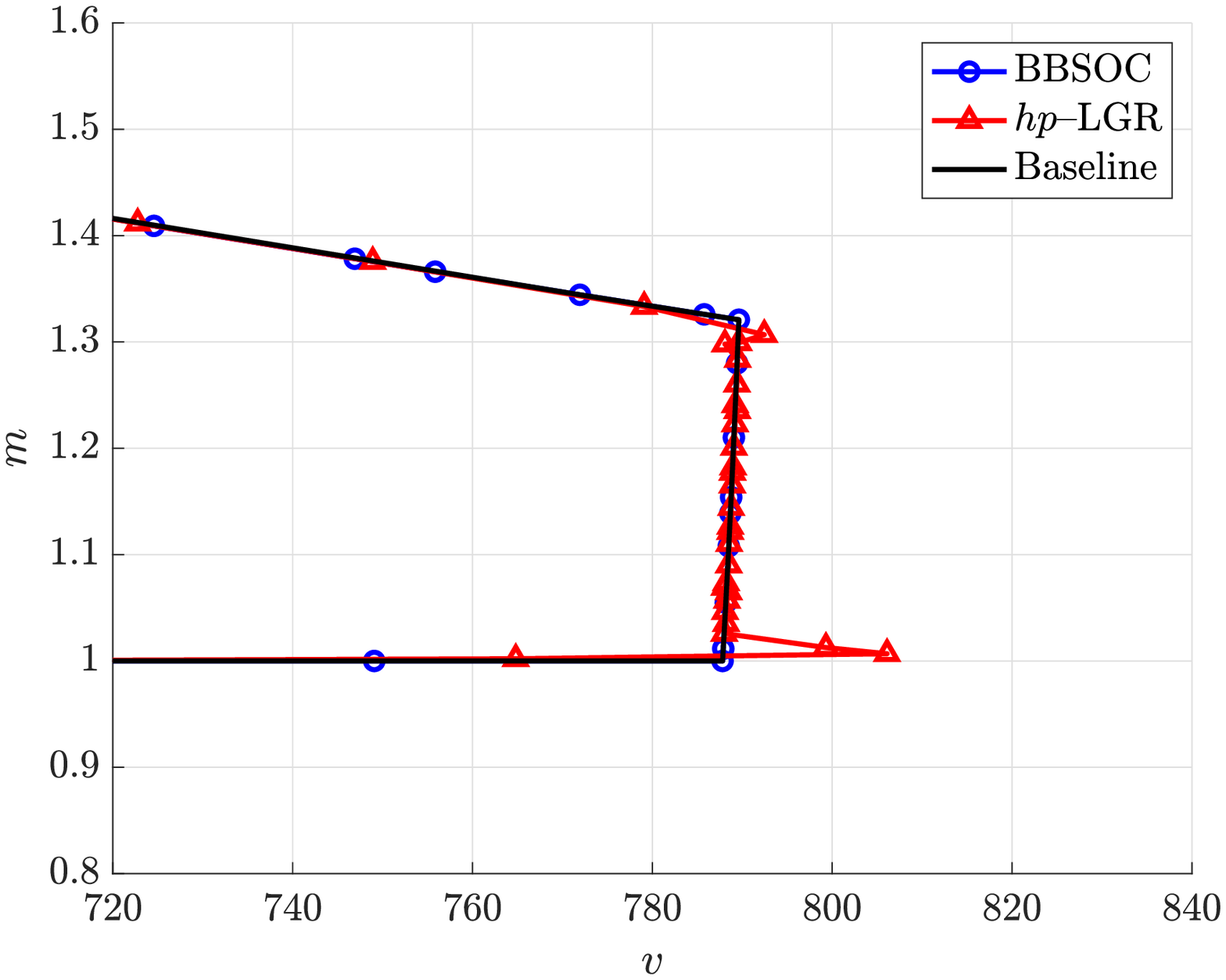}} \\
\caption{The singular surface for Example~$2$ is presented with the corresponding trajectory obtained using the BBSOC method, $hp$--LGR method, and the baseline solution. The numerical error accrued by the $hp$--LGR solution is shown in a zoomed in 2D plot of the trajectory in Fig.~\ref{fig:goddardRocket-SurfZoom}} 
\label{fig:goddardRocketSurf}
\end{figure}

It is well known that the singular surface in the optimal control problem of Eq.~\eqref{eq:goddardRocket} can be reduced to lie in the state space (that is, the singular surface is a three-dimensional surface in the space defined by $h$, $v$, and $m$).  In particular, the singular surface is defined as
\begin{equation}\label{eq:goddardRocket-surf}
	m = \frac{D_0v^2 \exp(-h/H)(1+v/c)}{g}.
\end{equation}
Figure~\ref{fig:goddardRocketSurf} shows the singular surface defined by Eq.~\eqref{eq:goddardRocket-surf} and the state solution obtained using the BBSOC method, $hp$--LGR method, and baseline solution.  It is observed in Fig.~\ref{fig:goddardRocket-SurfZoom} that the segment of the $hp$--LGR solution that corresponds to the singular interval does not lie in close proximity to the singular surface.  On the other hand, the portion of the trajectory corresponding to the singular interval using either the BBSOC or the baseline solution does lie on the singular surface.  

\begin{table}
\caption{\label{tab:goddardRocket} Comparison of computational results for Example~$2$.}
\centering
\begin{tabular}{ c " c c c c c c c c } 
\hline \hline
  &  $t_{s}^{[1]}$  & $t_{s}^{[2]}$  &  $t_f$  & $\C{J}^*$   & $\delta$    & $\epsilon$  &  $p$ & CPU~$[\textrm{s}]$  \\ \thickhline 
  BBSOC  & $13.751266$  & $21.987362$  & $42.887912$  & $-18550.87185$ & $7.86 \times 10^{-8}$ & $10^{-6}$ & $3$ & $2.86$ \\ 
  $hp$--LGR  & $13.541287$  & $22.213037$ & $42.887750$ & $-18550.87039$ & $-$  & $-$  & $-$ & $2.24$  \\
  Baseline  & $13.751270$  & $21.987363$ & $42.887912$ & $-18550.87186$ & $-$    &   $-$   & $-$ & $2.57$ \\
  \hline \hline
\end{tabular}
\end{table}

\subsection*{Example 3: Jacobson's Problem}

Consider the following problem posed in~\cite{Jacobson1970}:
\begin{equation}\label{eq:aly}
	\begin{array}{lcl}
    &\textrm{minimize}& \,\, \displaystyle \C{J} = \frac{1}{2	} \int_{0}^{t_f} (x_1^2(t) + x_2^2(t))\, dt, \vspace{0.5cm} \\
    &\textrm{subject to}& 
    \left\{
    \begin{array}{lclclclclcl}
		\dot{x_1}(t) &=& x_2(t)    &,&  x_1(0) &=& 0 &,& x_1(t_f) &=& \textrm{Free}, \\
		\dot{x_2}(t) &=& u(t)        &,&  x_2(0) &=& 1 &,& x_2(t_f) &=& \textrm{Free}, \\
		-1 &\leq& u(t) \,\,\,  \leq 1 &,&
	\end{array}
    \right.
    \end{array}
\end{equation} 
where $t_f = 5$.  The optimal control problem of Eq.~\eqref{eq:aly} has an analytic solution.  Consequently, the analytic solution can be used to assess the accuracy of the BBSOC method.  Further details on the derivation of the analytic solution to the example in Eq.~\eqref{eq:aly} can be found in~\cite{Aghaee2021}.  The singular control is $u_{\textrm{sing}}^*(t)=x_1(t),\; t\geq\approx 1.41376409$ and the analytic switch time $t_s^{[1]}=1.41376409$ is shown in Table~\ref{tab:aly}.

\begin{figure}
\centering
\vspace*{0.25cm}
\subfloat[\label{fig:aly-U}]{\includegraphics[scale=0.4]{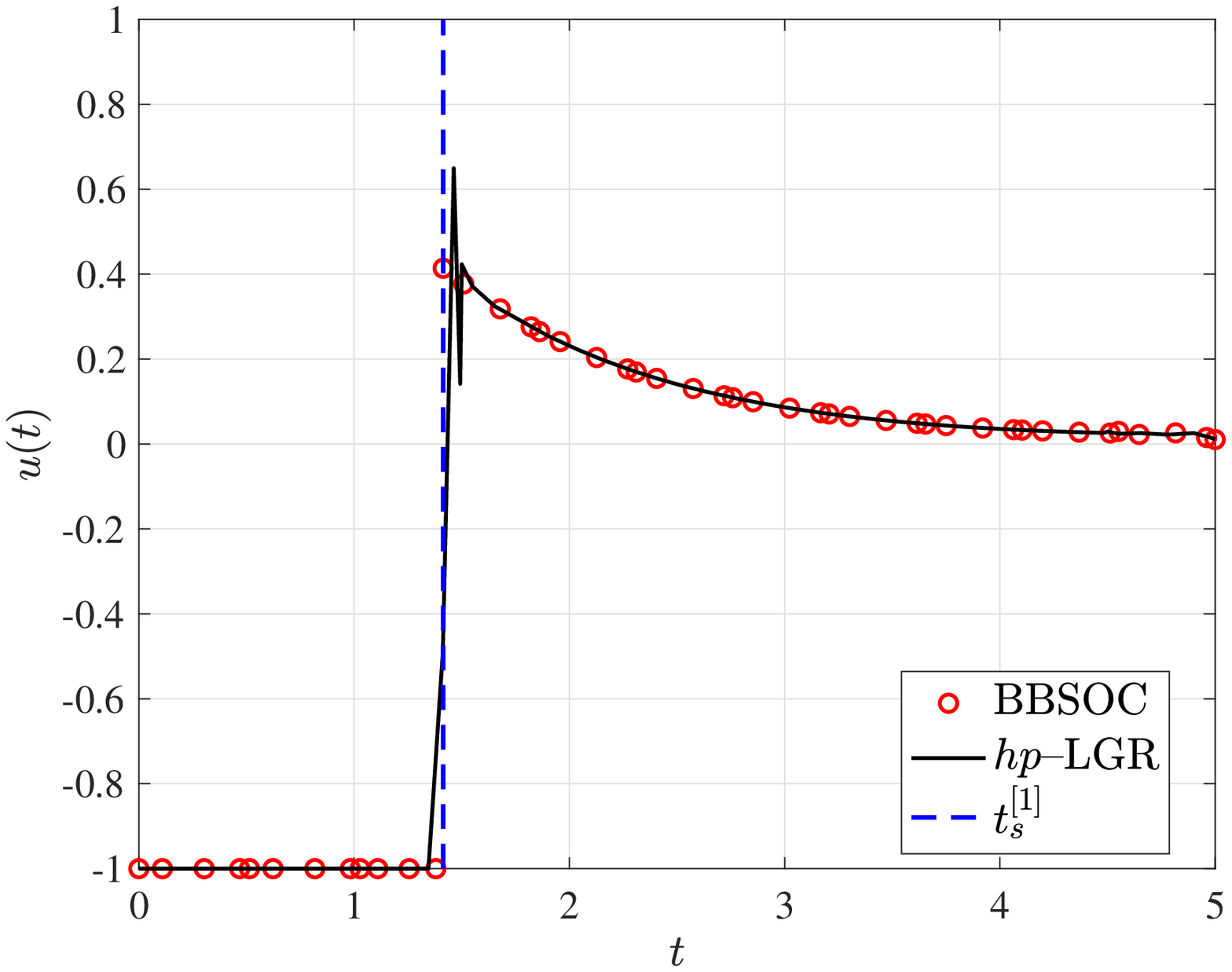}}~~\subfloat[\label{fig:aly-Hist}]{\includegraphics[scale=0.4]{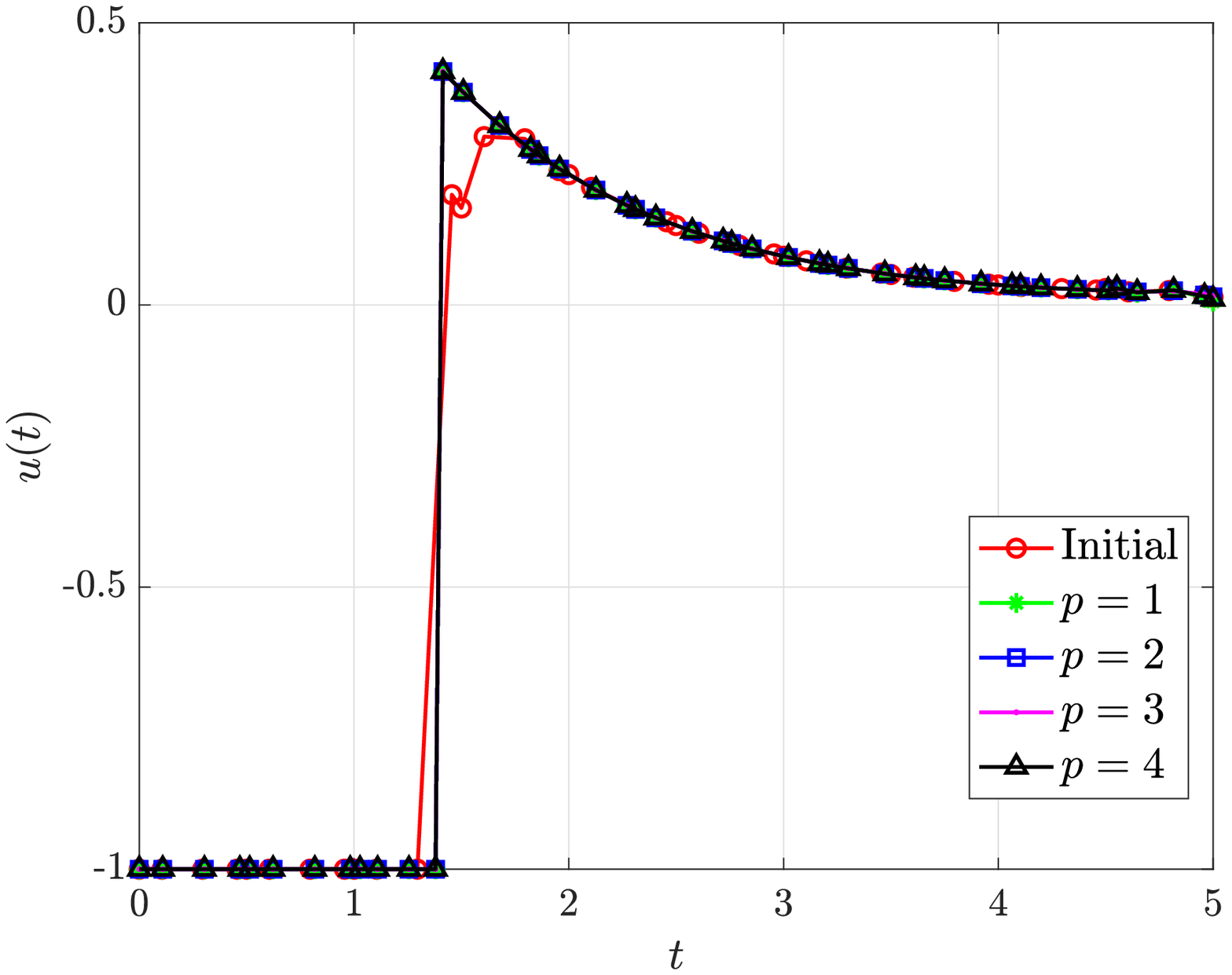}} \\ 
\caption{Control solution for Example~$3$ illustrating the numerical discrepancies between the solutions obtained by the BBSOC and $hp$--LGR methods. The control history of the regularization procedure shows the process of obtaining the singular control in Fig.~\ref{fig:aly-Hist}. } 
\label{fig:aly}
\end{figure}

Figure~\ref{fig:aly} shows the control solutions obtained when solving the optimal control problem of Eq.~\eqref{eq:aly} using both the BBSOC and $hp$--LGR methods.  Using the BBSOC method, it is seen that a high-accuracy approximation of both the switch time and the singular control are obtained after two iterations of the BBSOC method.  In contrast, the $hp$--LGR solution exhibits fluctuations in the neighborhood of the switch time.  Next, Fig.~\ref{fig:aly-Hist} shows the history of the control obtained on each iteration of the regularization method that is part of the BBSOC method.  As already indicated in Fig.~\ref{fig:aly} and further emphasized by Fig.~\ref{fig:aly-Hist}, these fluctuations present in the $hp$--LGR solution are eliminated because of the inclusion of the regularization term in the objective functional.  

Next, Table~\ref{tab:aly} provides a comparison of the switch time and objective values obtained using both the BBSOC and $hp$--LGR methods alongside the analytic solution (where it is noted that Table~\ref{tab:aly} also provides the parameter $\epsilon$ used in the regularization term of Eq.~\eqref{eq:penaltyTerm} along with the value $\delta$ obtained for the regularization term itself).   It is seen that the value of the analytic switch time is in close agreement with the switch time obtained from the BBSOC method.  Next, the computation times for both methods are also compared in Table~\ref{tab:aly}.  It is observed that the BBSOC method is more computationally efficient when compared with the $hp$--LGR method while simultaneously producing a more accurate solution.  

\begin{table}
\caption{\label{tab:aly} Comparison of computational results for Example~$3$}
\centering
\begin{tabular}{ c " c c c c c c c c } 
\hline \hline 
                &  $t_{s}^{[1]}$  & $\C{J}^*$   & $\delta$    & $\epsilon$  &  $p$ & CPU~$[\textrm{s}]$ \\ \thickhline 
BBSOC  & $1.41376404$   & $0.37699193$ & $6.31 \times 10^{-14}$ & $10^{-8}$ & $4$ & $2.02$  \\ 
$hp$--LGR  & $1.41194167$  & $0.37699503$ & $-$  & $-$  & $-$ & $2.13$  \\
Analytic  & $1.41376409$  & $-$ & $-$    &   $-$   & $-$  & $-$ \\
\hline \hline
\end{tabular}
\end{table}

\subsection*{Example 4:  Entry Vehicle Crossrange Maximization (Smooth Solution)}

Consider now the crossrange maximization of an entry vehicle problem given in Ref.~\cite{Betts2010} and whose formulation is taken from Ref.~\cite{Patterson2014}.  It is known that the solution to this optimal control problem is smooth.  For problems whose solutions are smooth, the BBSOC method should not identify any discontinuities or singular intervals.  Consequently, the method should not divide the entire domain into multiple domains nor should it apply regularization.  In this section the entry vehicle problem is solved using the BBSOC method.  While the details of this problem are omitted here for brevity, further details of the maximum crossrange entry vehicle problem can be found in Refs.~\cite{Betts2010,Patterson2014}.  

The results of applying the BBSOC method to the aforementioned maximum crossrange entry vehicle problem are provided in Table~\ref{tab:RLV}, where $M$ denotes the number of mesh refinement iterations and $K$ denotes the total number of collocation points.  In particular, it is seen from Table~\ref{tab:RLV} that the BBSOC method obtains the same solution as the optimal control software in Ref.~\cite{Patterson2014} while taking an average of $0.6$~s more CPU time to solve the problem.  In particular, the BBSOC method does not identify any discontinuities nor does it identify any intervals as singular.  The results obtained for this example demonstrate the ability of the BBSOC method to identify correctly the fact that the optimal control for the maximum crossrange entry vehicle problem is smooth and that the BBSOC method applies only static mesh refinement (see Section~\ref{sect:smooth}) as necessary.

\begin{table}[h]
\caption{\label{tab:RLV} Comparison of computational results for Example~$4$.}
\centering
\begin{tabular}{ c " c c c c } 
\hline \hline
                      			   &  $\C{J}^*$    & $M$  & $K$  & CPU~$[\textrm{s}]$   \\ \thickhline 
BBSOC       				   & $-0.5963$    & $4$   &  $105$ & $1.81$  \\ 
$hp$--LGR & $-0.5963$    & $4$  & $105$   & $1.17$    \\
\hline \hline
\end{tabular}
\end{table}

\section{Limitations of the BBSOC Method \label{sect:limitations}}

As with any computational method, the method of this paper has limitations.  The first limitation is that the approach depends upon the quality of the initial mesh that is supplied.  In particular, the solution on the initial mesh is used to ascertain the structure of the solution to the optimal control problem.  As a result, if the solution structure determined by the structure decomposition method described in Section \ref{sect:structureDecomp} is not representative of the optimal solution structure, the method may not produce a sufficiently accurate solution.  

Next, the regularization method described in Section \ref{sect:singular} is dependent upon the value of the regularization parameter $\epsilon$ as shown in Eq.~\eqref{eq:penaltyTerm}.  If $\epsilon$ is chosen to be too large, then the control obtained in a domain categorized as singular may not be a sufficiently accurate approximation of the singular control.  On the other hand, if $\epsilon$ is chosen to be too small, then the approximation of the singular control may be noisy (which would be as if no regularization was performed).  Next, if the structure detection method produces initial estimates of the switch times that do not lie within a reasonable proximity of the actual switch times or incorrectly identifies the actual switch times, then the structure obtained may be quite different from the optimal control structure. The detection of the control structure can be adjusted by changing the values of the parameters $\eta$ and $\mu$, as discussed in Sect.~\ref{sect:sructureDet}. While the default values of the parameters given in Sect.~\ref{sect:sructureDet} worked for the problems discussed in this paper, more complex problems may require more testing from a user standpoint.  

Finally, it is noted that the method of this paper is limited to problems where the control appears linearly in the Hamiltonian, that is, this paper does not consider problems where the solution is bang-bang and/or singular and the Hamiltonian is not linear in the control.  Furthermore, this method does not consider the inclusion of state or mixed state-input constraints. 


\section{Conclusions\label{sect:conc}}

A method has been described for solving bang-bang and singular optimal control problems using adaptive Legendre-Gauss-Radau (LGR) collocation.  First, the standard single-domain LGR collocation method is modified to be formulated as a multiple-domain method where the endpoints of each domain are variables in the optimal control problem.  Next, a structure detection method is developed that identifies intervals of the original domain as either regular, bang-bang, or singular.  This structure detection method is developed based on the sign of the switching function for problems where the control appears linearly in the optimal control Hamiltonian.  Based on the results of the structure detection, the original domain is partitioned into multiple-domains in a manner consistent with the aforementioned multiple-domain formulation.  For any domain that is categorized as regular, the control is set equal to either its lower or upper limit (based on the sign of the switching function).  For any domain that is categorized as singular (that is, the sign of the switching function in such a domain is indeterminate), a regularization method is employed to determine an accurate approximation of the singular control.  The various parts of the method, that is, the multiple-domain formulation, the structure detection, and the categorization of the various intervals, are combined into a unified method that solves bang-bang and singular optimal control problems.  It is also noted that, for problems where the solution contains no bang-bang or singular intervals, the method reverts to a standard single-domain LGR collocation method.  The method is demonstrated on four examples.  Three of the examples have either a bang-bang and/or singular solution while the fourth example contains no bang-bang or singular intervals.  It is found that the method efficiently produces accurate solutions to bang-bang and singular optimal control problems.

\section*{Data Availability}
The data that support the findings of this study are available from the corresponding author upon request.


\section*{Acknowledgments}\label{sect:ack}

The authors gratefully acknowledge support for this research from the U.S.~National Science Foundation under grants DMS-1819002 and CMMI-2031213, the U.S.~Office of Naval Research under grant N00014-19-1-2543, and from Lockheed-Martin Corporation under contract 4104177872.


\renewcommand{\baselinestretch}{1.0}
\normalsize\normalfont
\bibliographystyle{aiaa}


\end{document}